\documentclass{amsart}
\usepackage[T1]{fontenc}
\usepackage{textcomp}
\usepackage{amsfonts}
\usepackage{ae,aecompl}
\usepackage{amsmath}
\usepackage{amssymb}
\usepackage{amsthm}
\usepackage{epsf}

\theoremstyle{definition}

\newcommand{\ep}{\varepsilon}

\newcommand{\DtoN}{\Lambda_{\gamma}}

\newcommand{\bound}{\partial\Omega}

\providecommand{\abs}[1]{\lvert {#1} \rvert}
\providecommand{\norm}[1]{\lVert {#1} \rVert}
\newcommand{\diff}{\Omega_2\setminus\overline{\Omega_1}}
\newcommand{\pd}[2]{\frac{\partial {#1}}{\partial {#2}}}

\newcommand{\ra}{\rightarrow}

\newcommand{\T}{{\mathbf{t}}}
\newcommand{\Om}{\Omega}
\newcommand{\DOm}{\partial\Omega}

\newcommand{\R}{{\mathbb R}}

\newcommand{\C}{{\mathbb C}}
\renewcommand{\L}{{\mathcal L}}

\newcommand{\RR}{{\mathcal R}}

\title{Reconstructing conductivities \\with boundary corrected D-bar method}
\author{Samuli Siltanen and Janne P.~Tamminen}

\begin{document}

\begin{abstract}
\noindent
The aim of electrical impedance tomography is to form an image of the
conductivity distribution inside an unknown body using electric
boundary measurements. The computation of the image
from measurement data is a non-linear ill-posed inverse problem and
calls for a special regularized algorithm. One
such algorithm, the so-called D-bar method, is improved in this work
by introducing new computational steps that remove the so far 
necessary requirement that the conductivity should be constant near the boundary.
The numerical experiments presented suggest two conclusions. First, for most conductivities arising in medical imaging,
it seems the previous approach of using a best possible constant near the boundary is sufficient.
Second, for conductivities that have high contrast features at the boundary, the new approach produces 
reconstructions with smaller quantitative error and with better visual quality.
\end{abstract}

\maketitle


\section{Introduction}

\noindent
The aim of electrical impedance tomography ({\sc eit}) is to form an image of the
conductivity distribution inside an unknown body using electric
boundary measurements. Applications of {\sc eit} include medical imaging,
nondestructive testing and subsurface monitoring. See
\cite{cheney} for an overview of {\sc eit}. The computation of the image
from measurement data is a non-linear ill-posed inverse problem and
calls for a special regularized algorithm. In this work we improve one
such algorithm, the so-called D-bar method, by removing the so far
necessary requirement that the conductivity should be constant near
the boundary.

The mathematical model behind {\sc eit} is the {\em inverse conductivity
  problem} introduced by Calder\'on in \cite{calderonkys}. We discuss
here the inverse conductivity problem in the following two-dimensional form: let
$\Om_1=D(0,r_1)\subset\R^2$ be the disc with center at origin and radius
$r_1>0$ and consider a strictly positive, real-valued conductivity
$\sigma\in C^2(\overline{\Om_1})$. Maintaining a voltage distribution
$f$ on the boundary $\partial\Om_1$ creates a voltage potential $u$
that solves the following Dirichlet problem:
\begin{eqnarray}
  \left\{ \begin{array}{rcl}
  \nabla \cdot (\sigma \nabla u) &=& 0  \quad \rm{in} \quad \Om_1 \label{ody1},  \\
  u&=&f  \quad \rm{on} \quad \partial \Om_1.
\end{array} \right.
\end{eqnarray}
The resulting distribution of current through the boundary is
\begin{equation}
  \Lambda_\sigma f = \sigma\frac{\partial u}{\partial \nu}|_{\partial\Om_1}
\end{equation}
where  $\nu$ is the outward unit normal and $\Lambda_\sigma$ is the
Dirichlet-to-Neumann ({\sc dn}) map. Calder\'on's problem is
to reconstruct $\sigma$ from the knowledge of $\Lambda_\sigma$.

Many numerical methods have been suggested in the literature for the
reconstruction of $\sigma$ in the above setting. In this work we
concentrate on the so-called D-bar method. Other approaches include
linearization \cite{barber,JuttaJen,noser}, iterative regularization
\cite{dobson,dobson_santosa1,lechleiter,kindermannleitao}, statistical
inversion \cite{kaipio,RoyNichollsFox}, resistor network methods
\cite{BorceaDruskinVasquez}, convexification \cite{Klibanov}, layer
stripping \cite{eki,sylvester2D} and Teichm\"uller space methods
\cite{KLO1,KLO2,KLO3}. Also, there is a large body of work
concentrating on recovering partial information on $\sigma$, see
\cite{borcea,borcea2} for a survey.

Theoretical foundation of the D-bar method for {\sc eit} in dimension two
was introduced by Nachman in \cite{nachman}, where a constructive
proof for recovering $\sigma\in W^{2,p}(\Om_1)$ from $\Lambda_\sigma$
was given for $p>1$. The result was later sharpened by Brown and
Uhlmann \cite{brownuhlmann} to cover $W^{1,q}(\Om_1)$ conductivities
with $q>2$; the proof was augmented with constructive steps by Knudsen
and Tamasan in \cite{knudsenTamasan}. Finally, Astala and
P\"aiv\"arinta answered Calder\'on's question in its original form by
describing a constructive procedure for recovering $\sigma\in 
L^\infty(\Om_1)$ in  \cite{astala,astala2}; numerical implementation of this approach is described in \cite{AMPPS,AMPS}. Thus there are several
variants of the D-bar method for two-dimensional {\sc eit}. In dimension
three, the theory of D-bar reconstruction is developed in
\cite{towards3d,nachman88,novikov,sylvesteruhlmann,3Ddbar}. 

The above theoretical results on the D-bar method assume the knowledge of the infinite-precision data $\Lambda_\sigma$. 
However, the starting point of practical inversion is a finite-dimensional and
noisy approximation $L_{\sigma}^{\ep}$ to $\Lambda_\sigma$. Since the {\sc eit} problem is severely ill-posed,
or sensitive to measurement noise, any practical reconstruction method needs to be robust against errors in
 measurement data. The first robust D-bar algorithm (based on
 \cite{nachman}) was given in \cite{siltanen}, and it has been refined 
 and analysed in \cite{siam,FIST,SIAP,AIP2007}. The method has been
 successfully tested on a chest phantom in \cite{TMIdata} and on {\em
   in vivo} human chest data in \cite{ChestPaper}. Numerical D-bar
 method based on \cite{brownuhlmann,knudsenTamasan} was reported in
 \cite{knudsen}. The above methods are two-dimensional;
 three-dimensional computations are described in
 \cite{juttathesis,boverman}. Robustness is ensured in all of these
 practical D-bar reconstruction methods by truncating scattering data,
 a step that can be viewed as nonlinear low-pass filtering.

In addition to being robust against noise, a reliable 
{\sc eit} algorithm needs a regularization analysis. Such an analysis is provided
for the two-dimensional D-bar method in \cite{regularizedEIT}, where an explicit formula is given for choosing
the truncation radius as function of noise level. 

Practical D-bar methods
have been until now implemented by fitting an optimal constant to the
possibly nonconstant trace $\sigma|_{\bound_1}$. Quite good results
have been obtained both with laboratory data \cite{TMIdata} and {\em
  in vivo} patient data \cite{ChestPaper}. However, in applications
exhibiting large conductivity changes near or at the boundary the
constant-fitting approach may not be good enough. Our aim here is to
remove the assumption ``$\sigma\equiv 1$ near the boundary'' from the
two-dimensional D-bar algorithm based on \cite{nachman} using an
additional procedure we call \emph{boundary correction}.

Let us review the infinite-precision boundary correction procedure
given in \cite{nachman}. The starting point is the {\sc dn} map
$\Lambda_\sigma$ of a conductivity $\sigma\in W^{2,p}(\Om_1)$. Take
$r_2>r_1$ and set $\Om_2=D(0,r_2)$. The conductivity $\sigma$ is
extended outside $\Om_1$ by 
\begin{equation}\label{def:extension}
\gamma(x) = 
\begin{cases}  
\sigma(x), & \text{when }x\in\Om_1, \\
\tilde{\sigma}(x), & \text{when }x\in\diff, 
\end{cases}
\end{equation}
where we can choose any $\tilde{\sigma}\in
W^{2,p}(\Om_2\setminus\overline{\Om_1})$ with the properties
$\tilde{\sigma}|_{\bound_1}=\sigma|_{\bound_1}$ and
$(\partial\sigma/\partial\nu)|_{\partial\Om_1} =
(\partial\tilde{\sigma}/\partial\nu)|_{\partial\Om_1}$ and
$\tilde{\sigma}\equiv1$ near $\bound_2$. This way $\gamma\in
W^{2,p}(\Om_2)$ whenever $\sigma\in W^{2,p}(\Om_1)$.
Define two Dirichlet problems:
\begin{equation}
\left\{
\begin{array}{rcl}
\nabla\cdot(\tilde{\sigma}\nabla u_j) &=& 0 \quad\rm{in}\quad\diff,\quad j=1,2 \\
u_j &=& f_j\quad\rm{on}\quad\bound_j \label{extensionproblem}\\
u_j &=& 0 \quad\rm{on}\quad\bound_i, \quad i=1,2, \quad i\neq j.
\end{array}\right .
\end{equation}
Four new {\sc dn} maps in $\diff$ can be defined by
\begin{equation}\label{fourDNs}
\Lambda^{ij}f_j = \tilde{\sigma}\pd{u_j}{\nu}|_{\bound_i},\quad i,j=1,2.
\end{equation}
By proposition 6.1 of \cite{nachman} we can use (\ref{fourDNs}) to write $\DtoN$ in terms of $\Lambda_\sigma$:
\begin{equation}\label{LgintermsofLs}
\DtoN = \Lambda^{22}+\Lambda^{21}(\Lambda_\sigma-\Lambda^{11})^{-1}\Lambda^{12}.
\end{equation}
The boundary corrected D-bar method for $\sigma\in W^{2,p}(\Om_1)$,
assuming infinite-precision data, is described as follows in
\cite[Section 6]{nachman}:
\begin{tabbing}
\\
  $\quad$ \=(a) \= {\bf Reconstruction at the boundary.}  Recover the trace $\sigma|_{\partial\Om_1}$ and \\ \>\>
  the normal derivative $(\partial\sigma/\partial\nu)|_{\partial\Om_1}$ from $\Lambda_\sigma$.
  \\\\
  \>(b)\>{\bf Extension of conductivity.} 
  Using (a) and (\ref{def:extension}), extend the conductivity\\\>\> to $\gamma \in W^{2,p}(\Om_2)$
  such that  $\gamma\equiv 1$ near $\partial\Om_2$ and $\inf_{x\in\Om_2} \gamma(x)\geq c>0$;
  \\\\
  \>(c) \>{\bf Calculation of outer {\sc dn} map.} 
  Write the {\sc dn} map $\DtoN$ of $\gamma \in W^{2,p}(\Om_2)$ \\\>\> in terms
  of $\Lambda_\sigma$ using (\ref{LgintermsofLs});
  \\ \\
  \>(d) \>{\bf Reconstruction using the D-bar method.} 
  Reconstruct $\gamma \in W^{2,p}(\Om_2)$ \\\>\> from the infinite-precision data $\DtoN$ following \cite{nachman}.
  \\
\end{tabbing}

The practical starting point of reconstruction is the noisy
approximate data $L_\sigma^\ep$, and the Steps (a--d) above cannot be
directly followed. We suggest the following robust procedure for
boundary correction:
\begin{tabbing}
\\
  $\quad$ \=(a$^\prime$) \=  {\bf Approximate reconstruction at the boundary.} Recover
  numerically \\  \>\> a smooth function $g\in C^\infty(\partial\Om_1)$ with
  the property $g\approx\sigma|_{\partial\Om_1}$
  as explained \\\>\> in \cite{boundary1}. Omit recovering $(\partial\sigma/\partial\nu)|_{\partial\Om_1}$
since it is an unstable step \cite{boundary1};
  \\  \\
  \>(b$^\prime$) {\bf Simple extension of conductivity.}  Construct  $\tilde{\sigma}\in
  C^2(\overline{\Om_2\setminus\Om_1})$ \\\>\> satisfying
  $\inf_{x\in\Om_2\setminus\overline{\Om_1}} \tilde{\sigma}(x)\geq
  c>0$
   and $\tilde{\sigma}|_{\partial\Om_1}=g$ and $\tilde{\sigma}\equiv
   1$ near $\partial\Om_2$. \\\>\>Use $\tilde{\sigma}$ in (\ref{def:extension}) to extend the
   conductivity to $\gamma\in L^\infty(\Om_2)$;
  \\  \\
  \>(c$^\prime$) {\bf Approximate calculation of outer {\sc dn} map.}
  Write approximate \\\>\> {\sc dn} map $L_\gamma$ in terms of $L_\sigma$ using a matrix approximation to (\ref{LgintermsofLs});
  \\  \\
  \>(d$^\prime$) {\bf Reconstruction using regularized D-bar method.} Reconstruct $\gamma$ \\\>\> from
  $L_\gamma$ using the regularized D-bar method described in \cite{regularizedEIT}.\\
\end{tabbing}
The main concern about the procedure (a$^\prime$--d$^\prime$)  is that
the extension of $\sigma$ to $\gamma$ will be in general discontinuous
at $\bound_1$, and thus $\gamma$ violates the assumptions of the D-bar
method used in (d$^\prime$). However, there is both theoretical and
experimental evidence suggesting that the step (d$^\prime$) should
give reasonable results even in this case \cite{TMIdata,ChestPaper,AIP2007,SIAP}. Another potential problem
arises from the inverse operator in formula (\ref{LgintermsofLs}), as the proof of invertibility \cite[Proposition 6.1]{nachman} in the extended conductivity produced by step (b$^\prime$). One possibility would be to use \cite[Lemma 2.1.3.]{KnudsenPhD} instead of (\ref{LgintermsofLs}) as the basis of step (b$^\prime$). However, in our computational experiments the use of (\ref{LgintermsofLs}) seems not to be a problem even in the case of discontinuous conductivity extensions.

This paper should be viewed as a report of computational experiments suggesting the practical usefulness of the boundary correction step in applications where the conductivity varies strongly near the boundary. Hopefully the computational results presented below will act as motivation for further theoretical study of practical imaging algorithms for {\sc eit}. 

We remark that the boundary correction method is applied in this paper only in the case of $\Omega_1$ being a disc. This
is not a serious lack of generality, though: we presume that other domains than discs could be treated combining the
methods described in \cite{murphy2007,MMN2007} with the boundary correction.

This paper is organized as follows. We present our method of
simulating continuum model {\sc eit} data in Section \ref{sec:Simulation}. The
details of implementation of Steps (a$^\prime$) and (b$^\prime$) and
(c$^\prime$) are discussed in Sections \ref{sec:boundary} and
\ref{sec:extension} and \ref{sec:outerDN}, respectively. A brief
outline of the regularized D-bar method is given in Section
\ref{sec:regDbar}. Our practical boundary correction method is
illustrated by numerical examples in Section \ref{sec:numerical}, and
we conclude our results in Section \ref{sec:conclusion}.

\section{Simulation of measurement data}\label{sec:Simulation}

\noindent
Let $\RR_\sigma: \widetilde{H}^{-1/2}(\partial\Om_1)\ra
\widetilde{H}^{1/2} (\partial\Om_1)$ denote the Neumann-to-Dirichlet map
of  $\sigma$, where $\widetilde{H}^s$ spaces consist of $H^s$
functions with mean value zero. We have $\RR_\sigma g = u|_{\partial\Om_1}$, where
$u$ is the unique $H^1(\Om_1)$ solution of the Neumann problem 
$$
  \nabla\cdot\sigma\nabla u = 0\mbox{ in }\Om_1,\qquad
  \gamma\frac{\partial u}{\partial \nu} = g \mbox{ on }\partial\Om_1,
$$
satisfying $\int_{\partial\Om_1}u ds=0$. We note two key
equalities concerning $\Lambda_\sigma$ and $\RR_\sigma$. Define a
projection operator $P\phi:=|\partial\Om_1|^{-1}\int_{\partial\Om_1}\phi$. Then for any $f\in
H^{1/2}(\partial\Om_1)$ we have $P \Lambda_\sigma f =
|\partial\Om_1|^{-1}\int_{\partial\Om_1}  \sigma\frac{\partial
  u}{\partial\nu} = 
\int_{\Om_1} \nabla\cdot  \sigma\nabla u = 0,$ so actually
$\Lambda_\sigma:H^{1/2}(\partial\Om_1) \ra \widetilde{H}^{-1/2}
(\partial\Om_1)$. From the definitions of $\Lambda_\sigma$ and 
$\RR_\sigma$ we now have 
\begin{eqnarray}
  \label{idenLR}
  \Lambda_\sigma \RR_\sigma &=& I \mbox{\ \ \ \ \ \ } \qquad:
  \widetilde{H}^{-1/2}(\partial\Om_1)\ra \widetilde{H}^{-1/2}
  (\partial\Om_1), 
  \\
  \label{idenRL}
  \RR_\sigma \Lambda_\sigma  &=& I - P \qquad:
  H^{1/2}(\partial\Om_1)\ra \widetilde{H}^{1/2} (\partial\Om_1).
\end{eqnarray}

Given $\sigma$ and $N>0$, we define a matrix
$R_\sigma:\C^{2N}\ra\C^{2N}$ as follows. We use a truncated
orthonormal trigonometric basis for representing functions defined at the
boundary $\bound_j$:
\begin{equation}\label{basisfunctions}
\phi^{(n)}_j(\theta) = \frac{1}{\sqrt{2\pi r_j}}e^{in\theta},\quad
n=-N,...,N,\quad j=1,2.
\end{equation}

Note that $\int_{\bound_j}\phi^{(n)}_jds =0$ for $n\neq0$. Then solve the Neumann problem
\begin{equation}\label{neumann}
  \nabla\cdot\sigma\nabla u^{(n)}_1 = 0\mbox{ in }\Om_1,\qquad
  \sigma\frac{\partial u^{(n)}_1}{\partial \nu} = \phi^{(n)}_1 \mbox{ on }\partial\Om_1,
\end{equation}
with the constraint $\int_{\bound_1} u^{(n)}_1ds =0.$ Define
$R_\sigma=[\widehat{u}(\ell, n)]$ by
\begin{equation}
   \widehat{u}(\ell, n) = \int_{\partial\Om_1} u^{(n)}_1 \overline{\phi^{(\ell)}_1} ds.
\end{equation}
Here $\ell$ is the row index and $n$ is the column index.

The matrix $R_\sigma$ represents the operator $\RR_\sigma$
approximately. We add simulated measurement noise by defining
\begin{equation}\label{noiselevel}
  R_\sigma^\ep := R_\sigma + cE,
\end{equation}
where ${E}$ is a $2N\times 2N$ matrix with random entries
independently distributed according to the Gaussian normal density
$\mathcal{N}(0,1)$. The constant $c>0$ is adjusted so that
$\|R_1^\ep-R_1\|/\norm{R_1}$, where $\norm{\cdot}$ is the standard matrix norm and $R_1$ is the ND -map for the unit conductivity,
is greater than the relative error caused by FEM and of the same order of magnitude as
0.0017\% (signal to noise -ratio of 95.5 dB) , the noise level of the {\sc ACT}3 impedance tomography imager
of Rensselaer Polytechnic Institute \cite{ACTIII}.

We can now easily compute the corresponding noisy matrix
representation $L_\sigma^\ep$ for the {\sc dn} map $\Lambda_\sigma$. Namely,
define 
$$
  \widetilde{L_\sigma^\ep} := (R_\sigma^\ep)^{-1};
$$
then $\widetilde{L_\sigma^\ep}$ is a matrix of size $2N\times 2N$. We should
add appropriate mapping properties for constant basis functions at the
boundary according to the facts 
$$
  \Lambda_\sigma 1=0, \qquad \int_{\partial\Om_1}\Lambda_\sigma f ds =0.
$$
This is achieved simply by setting (in Matlab notation)
\begin{equation}\label{zero_blocks}
  L_\sigma^\ep := \left[\begin{array}{ccc}
  \widetilde{L_\sigma^\ep}(1\!\!:\!\!N,1\!\!:\!\!N) & 0 & \widetilde{L_\sigma^\ep}(1\!\!:\!\!N,(N+1)\!\!:\!\mbox{end})\\
  \\
  0 & 0 & 0\\
  \\
  \widetilde{L_\sigma^\ep}((N+1)\!\!:\!\mbox{end},1\!\!:\!\!N) &\quad 0\quad &
  \widetilde{L_\sigma^\ep}((N+1)\!\!:\!\mbox{end},(N+1)\!\!:\!\mbox{end}) 
  \end{array}\right],
\end{equation}
where the zero block matrices above have various (but obvious) sizes.

\section{Approximate reconstruction at the boundary}\label{sec:boundary}

\noindent
The trace $\sigma|_{\bound_1}$ can be approximately
reconstructed in the following way \cite{boundary1}. Define $h_{M,\beta}(\theta) =
e^{iM\theta}\eta(\theta-\beta)$, where 

\begin{equation}
  \label{cutoff_function}
  \eta(\theta) =
  \begin{cases}
    d(\kappa\theta-\pi/2)^\alpha(\kappa\theta+\pi/2)^\alpha\cos(\kappa\theta),
    & \textrm{for } -\pi/(2\kappa)<\theta<\pi/(2\kappa),\\ 
    0, & \textrm{otherwise}
  \end{cases}
\end{equation}
is a non-negative cut-off function satisfying
$\int_{\bound_1}\eta^2(\theta)d\theta=1$. Now the mollified trace 
 $(\sigma\eta^2)|_{\bound_1}(\beta)$ can be calculated with
\begin{equation}
  \label{boundary_equation}
  \int_{\bound_1}\sigma\eta^2ds =
  \lim_{M\rightarrow\infty}\frac{1}{M}\int_{\bound_1}\overline{h_{M,\beta}}\Lambda_\sigma h_{M,\beta}ds.
\end{equation}
We get the approximation $g\approx\sigma|_{\bound_1}$ by
calculating \eqref{boundary_equation} with different angles $\beta$
and using a finite $M$ in the right side of \eqref{boundary_equation}.

Another approach to reconstructing $\sigma|_{\bound_1}$ is the layer stripping method introduced in \cite{eki}.

\section{Simple extension of conductivity}\label{sec:extension}

\noindent
The starting point here is a given approximation $g:\bound_1\ra\R$ to
the trace $\sigma|_{\bound_1}$ of the conductivity $\sigma\in
C^2(\overline{\Om_1})$ of interest. The aim is to construct a strictly
positive conductivity $\tilde{\sigma}:\Om_2\setminus\overline{\Om_1}$
satisfying $\tilde{\sigma}|_{\bound_1}=g$ and $\tilde{\sigma}\equiv 1$
near the outer boundary $\bound_2$, and then use formula
(\ref{def:extension}) to define $\gamma$.

We extend $\sigma$ to $\gamma$ using the following extension in polar
coordinates:
\begin{align}\label{sigma_tilde}
  \gamma(\rho,\theta) =
  \begin{cases}
     \sigma(\rho,\theta),& \rho \leq r_1,\\
      (g(\theta)-1)f_m(\rho)+1,& r_1< \rho \leq r_e,\\
     1,&  r_e < \rho \leq r_2,\\
  \end{cases}
\end{align}
where $r_1<r_e<r_2$ is some radius and $f_m(\rho)\geq0$ is a suitable third-degree polynomial satisfying
$f_m(r_1) = 1$ and $f_m(r_e)=0$. Note that $\gamma$ is twice
continuously differentiable apart from possible discontinuity at
$\rho=r_1$, and equals constant $1$ in the annulus  $r_e < \rho <
r_2$.

\section{Approximate calculation of outer {\sc dn} map}\label{sec:outerDN}

\noindent
Using the functions \eqref{basisfunctions}, a given function
$f:\bound_i\ra \C$ can be approximately represented by the vector
$$
 \vec{f} =
 [\hat{f}(-N),\hat{f}(-N+1),\dots,\hat{f}(N-1),\hat{f}(N)]^T,
\qquad \hat{f}(n)=\int_{\partial\Om_i} f\overline{\phi^{(n)}_i}ds,
$$
and the {\sc dn} maps $\Lambda^{ij}$ can be approximated by the
matrices $L^{ij}=[\widehat{g}_{ij}(\ell,n)]$ with 
\begin{equation}
  \widehat{g}_{ij}(\ell,n) =
  \int_{\bound_j}\tilde{\sigma}\pd{u^{(n)}_j}{\nu}|_{\bound_i}\overline{\phi^{(\ell)}_j} dS,
\end{equation}
where $u^{(n)}_j$ denotes the solution to \eqref{extensionproblem} with
$u^{(n)}_j|_{\bound_j} = \phi^{(n)}_j$. Again $\ell$ is the row index and
$n$ is the column index. Now the matrix $L_\gamma^\ep$ can be
calculated by
\begin{equation}
  \label{L_gamma}
  L_\gamma^\ep = L^{22}+L^{21}(L_\sigma^\ep-L^{11})^{-1}L^{12},
\end{equation}
provided that the matrix $L_\sigma^\ep-L^{11}$ is invertible.  Formula (\ref{L_gamma}) is a finite-dimensional
approximation to (\ref{LgintermsofLs}).

\section{Regularized D-bar method}\label{sec:regDbar}

\noindent In this section we explain how to reconstruct a conductivity
$\gamma$ in a regularized way from a noisy measurement matrix
$L^\ep_\gamma$ under the assumptions $\gamma\in C^2(\overline{\Om_2})$ and
$\gamma\equiv 1$ in a neighborhood of $\partial\Om_2$.

If we had the infinite-precision data $\Lambda_\gamma$ at our
disposal, we could follow the reconstruction procedure in
\cite{nachman}. 
First we would solve the boundary integral equation
\begin{align}  \label{eq:psieps}
  \psi(\,\cdot\,,k)|_{\DOm_2} = e^{ikx} - S_k(\Lambda_\gamma-\Lambda_1)\psi(\,\cdot\,,k)|_{\DOm_2},
\end{align}
in the Sobolev space $H^{1/2}(\partial\Om_2)$ for all
$k\in\C\setminus\{0\}$. In formula (\ref{eq:psieps}), $S_k$ is a
single-layer operator 
$$
  (S_k \phi)(x) :=\int_{\partial\Om_2}G_k(x-y)\phi(y)ds,
$$
where $G_k$ is Faddeev's Green function defined by 
$$
  G_k(x) := e^{ikx}g_k(x),\quad g_k(x) := \frac{1}{(2\pi)^2}
  \int_{\R^2}\frac{e^{ix\cdot\xi}}{|\xi|^2+2k(\xi_1+i\xi_2)}d\xi.
$$
Once equation (\ref{eq:psieps}) had been solved, we would substitute the result into 
\begin{align}\label{Tdef}
  \T(k)=\int_{\partial\Om_2}e^{i\bar{k}\bar{x}} (\Lambda_\gamma-\Lambda_1)\psi(\,\cdot\,,k)ds,
\end{align}
where $\T$ is called the scattering transform, and $\Lambda_1$ is the
{\sc dn} map for the unit conductivity. For each fixed $x\in
\Om$, we would solve the following integral formulation of the D-bar
equation: 
    \begin{align}\label{eq:int}
      \mu(x,k)=  1+\frac{1}{(2\pi)^2}\int_{\R^2}
      \frac{\T(k^\prime)}{(k-k^\prime)\bar{k}^\prime} e^{i(k^\prime x+\overline{k^\prime}\overline{x})}
      \overline{\mu(x,k^\prime)}dk^\prime_1 dk^\prime_2;
    \end{align}
then the conductivity would be perfectly reconstructed as $\gamma(x)=\mu(x,0)^2$.

However, since our starting point in practice is the matrix
$L_\gamma^\ep$, we need to regularize the above ideal approach as
explained in \cite{regularizedEIT}. The matrices $L_\gamma^\ep$ and $L_1$
we already have, and a matrix representation $\mathbf{S}_k$ for the
single-layer operator $S_k$ can be computed numerically, provided we have numerical evaluation routines for $g_k(x)$, see \cite{IS2004}. We expand
$e^{ikx}|_{\partial\Om_2}$ as a vector $\vec{g}$ in our
finite trigonometric basis (\ref{basisfunctions}) and set 
\begin{equation}\label{def:psiBIE}
  \vec{\psi}_k := [I+\mathbf{S}_k(L_\gamma^\ep-L_1)]^{-1}\vec{g}.
\end{equation}
for $k$ ranging in a fine grid inside the disc $|k|<R$, where the
truncation radius $R>0$ is ideally chosen according to the size of
noise. The choice of $R$ falls outside the scope of this paper, so we
will compute below reconstructions with $R$ ranging in an interval.
We define the truncated scattering transform by 
\begin{align}\label{eq:Teps}
    \T_R(k)=
  \begin{cases}
    \int_{\partial\Om_2}e^{i\bar{k}\bar{x}}
    \mathcal{F}^{-1}((L^\ep_\gamma-L_1)\vec{\psi}_k)(x) ds& \mbox{for }|k|<R, \\ 0, &
    \text{otherwise},
  \end{cases}
\end{align}
where $\mathcal{F}^{-1}$ denotes transforming from the Fourier series
domain to the function domain. Finally we use the numerical algorithm
in \cite{FIST} to solve equation (\ref{eq:int}) with $\T$ replaced by
$\T_R$ and denote the solution by $\mu_R(x,k)$. Then $\gamma(x)\approx
\mu_R(x,0)^2$.

\section{Numerical results}\label{sec:numerical}

\noindent  
We define several conductivity distributions $\sigma\in L^\infty(\Om_1)$ on the
unit disc $\Om_1=D(0,r_1)=D(0,1)$  and compare reconstructions computed
with and without the boundary correction procedure.

Before proceeding with the examples, though, we need to choose an optimal radius $r_2$ to be used in the boundary correction step. We do this by examining numerically the simple case of the unit conductivity $\sigma\equiv 1$ and using the procedure (a$^\prime$-c$^\prime$) explained in the introduction. The numerical parameters used in this procedure are the same as in the example reconstructions, and they are given later in this chapter.

We take  $N=16$ and simulate non-noisy ND map $R_1$ using the finite element method with 1048576 triangles in $\Om_1$ as explained in Section \ref{sec:Simulation}. Using the standard square norm for matrices, this yields $\epsilon_{\textrm{fem}}=\norm{R_1^{\textrm{th}}-R_1}/\norm{R_1^{\textrm{th}}}\approx0.0000173$,
where $R_1^{\textrm{th}}$ is the analytically calculated ND matrix for the unit
conductivity. Furthermore, we construct noisy ND map $R^\epsilon_1$ with formula (\ref{noiselevel}) and $c=0.00001$, giving\\ $\|R^\epsilon_1-R_1\|/\norm{R_1}\approx 0.0001>\epsilon_{\textrm{fem}}$.

\begin{figure}[p]
\begin{picture}(110,60)
\epsfxsize=5.5cm
\put(0,0){\epsffile{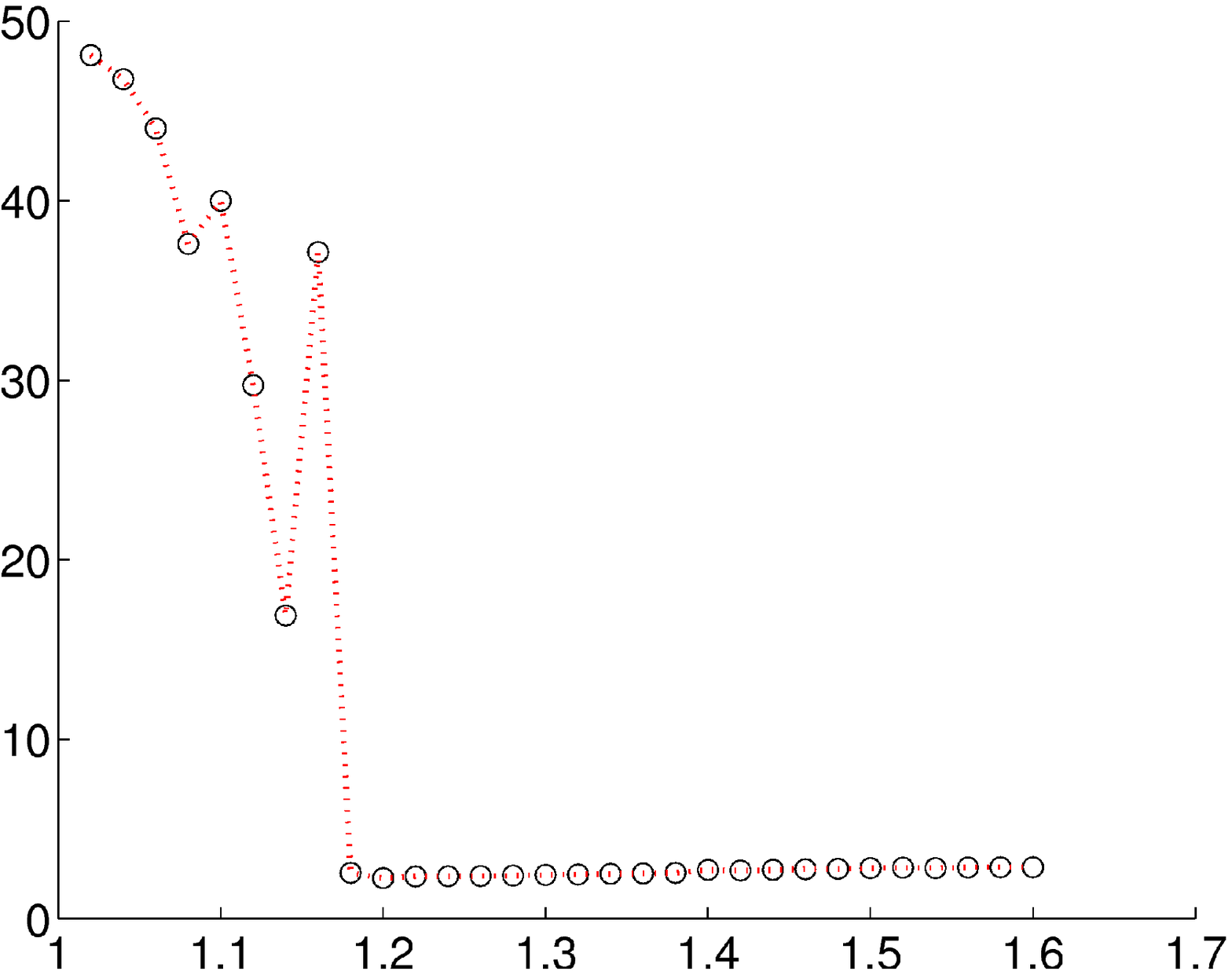}}
\epsfxsize=5.5cm
\put(55,0){\epsffile{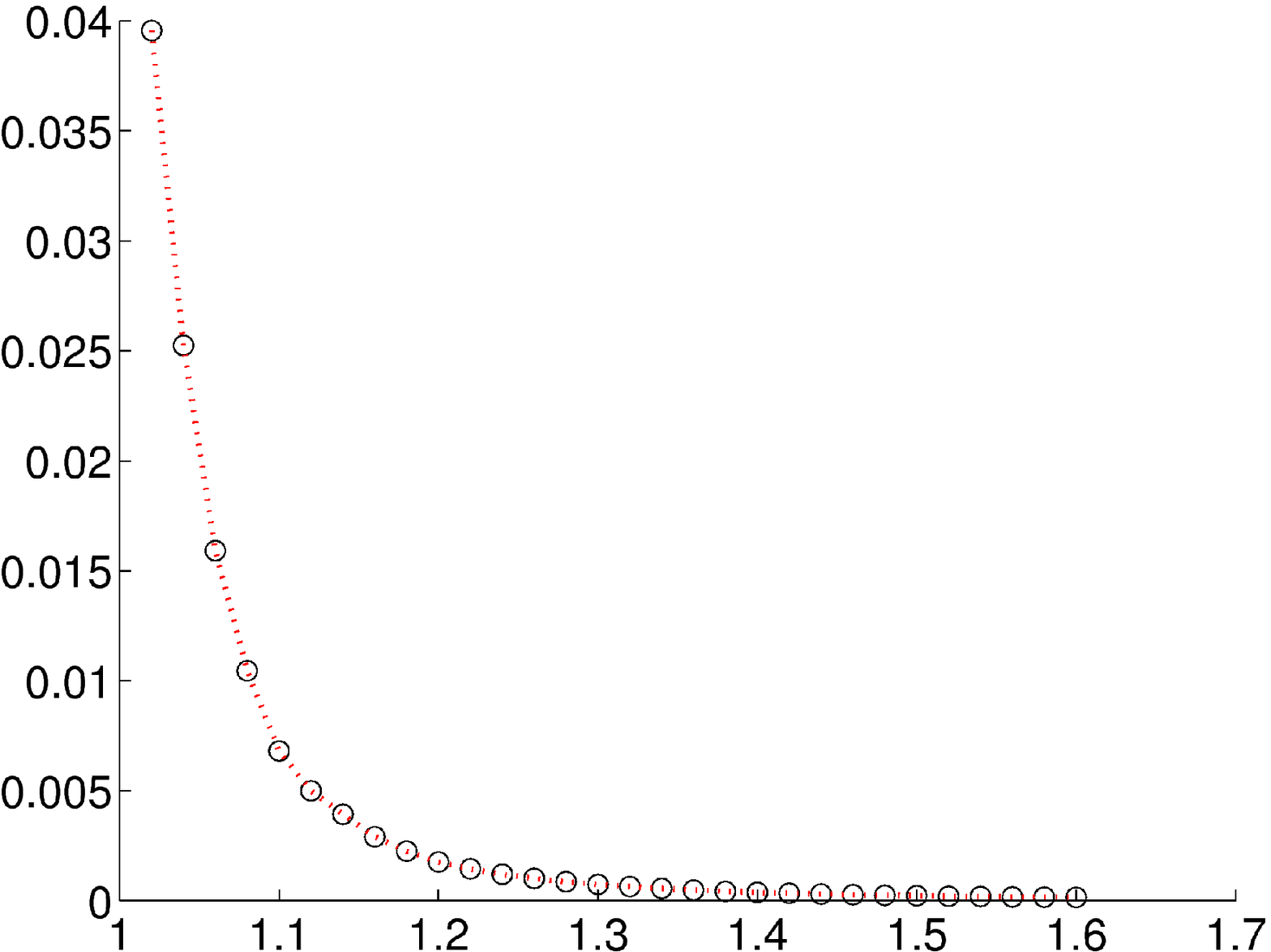}}
\put(110,-2){\tiny$r_2$}
\put(-3,40){\tiny \%}
\put(2,48){$\norm{L^\epsilon_{\gamma=1}-L^2_{\gamma=1}}/\norm{L^2_{\gamma=1}}$}
\put(60,48){$\norm{L^\epsilon_{\gamma=1}-L_{\gamma=1}}/\norm{L_{\gamma=1}}$}
\end{picture}
\caption{\label{fig:r2_test} Left: relative error $\norm{L^\epsilon_{\gamma=1}-L^2_{\gamma=1}}/\norm{L^2_{\gamma=1}}$ as a function of $r_2$. Here $L^\epsilon_{\gamma=1}$ is the {\sc dn} map on $\bound_2$ computed from noisy ND map using formula  (\ref{L_gamma}) and $L^2_{\gamma=1}$ is the {\sc dn} map calculated directly on $\bound_2$. Here $\|\,\cdot\,\|$ denotes the standard square norm for matrices. By $\gamma=1$ we mean $\sigma=1$ extended by \eqref{sigma_tilde}. Right: relative error $\norm{L^\epsilon_{\gamma=1}-L_{\gamma=1}}/\norm{L_{\gamma=1}}$ as a function of $r_2$. Here $L_{\gamma=1}$ is the {\sc dn} map computed from non-noisy ND map using formula  (\ref{L_gamma}).}
\end{figure}

To avoid notational clashes, we denote by
\begin{itemize}
\item[$L^\epsilon_{\gamma=1}$] the {\sc dn} map on $\bound_2$ computed from noisy ND map using formula  (\ref{L_gamma}),
\item[$L_{\gamma=1}$] the {\sc dn} map computed from non-noisy ND map using formula  (\ref{L_gamma}), 
\item[$L^2_{\gamma=1}$] the {\sc dn} map computed directly on $\bound_2$,
\end{itemize}
where by $\gamma=1$ we mean the conductivity $\sigma=1$ extended by \eqref{sigma_tilde}. The left plot in Figure \ref{fig:r2_test} shows the behaviour of the error $\norm{L^\epsilon_{\gamma=1}-L^2_{\gamma=1}}/\norm{L^2_{\gamma=1}}$  as a function of $r_2$. The condition number of the matrix $L^\epsilon_{\sigma=1}-L^{11}$ ranges between 1 and 20.  It seems that we should choose $r_2\geq 1.2$. Further, the right plot in Figure \ref{fig:r2_test} shows the behaviour of the error 
$\norm{L^\epsilon_{\gamma=1}-L_{\gamma=1}}/\norm{L_{\gamma=1}}$   as function of $r_2$. The error decreases as $r_2$ grows; it shows how the data measured on $\bound_1$ contributes less and less to $L^\epsilon_{\gamma=1}$ as $r_2$ gets larger. This observation is in agreement with the known fact that in {\sc eit} it is more difficult to obtain information from the deeper parts of the object \cite{distinguishability}.

\begin{figure}[p]
\begin{picture}(120,190)
\epsfxsize=4.5cm
\put(0,0){\epsffile{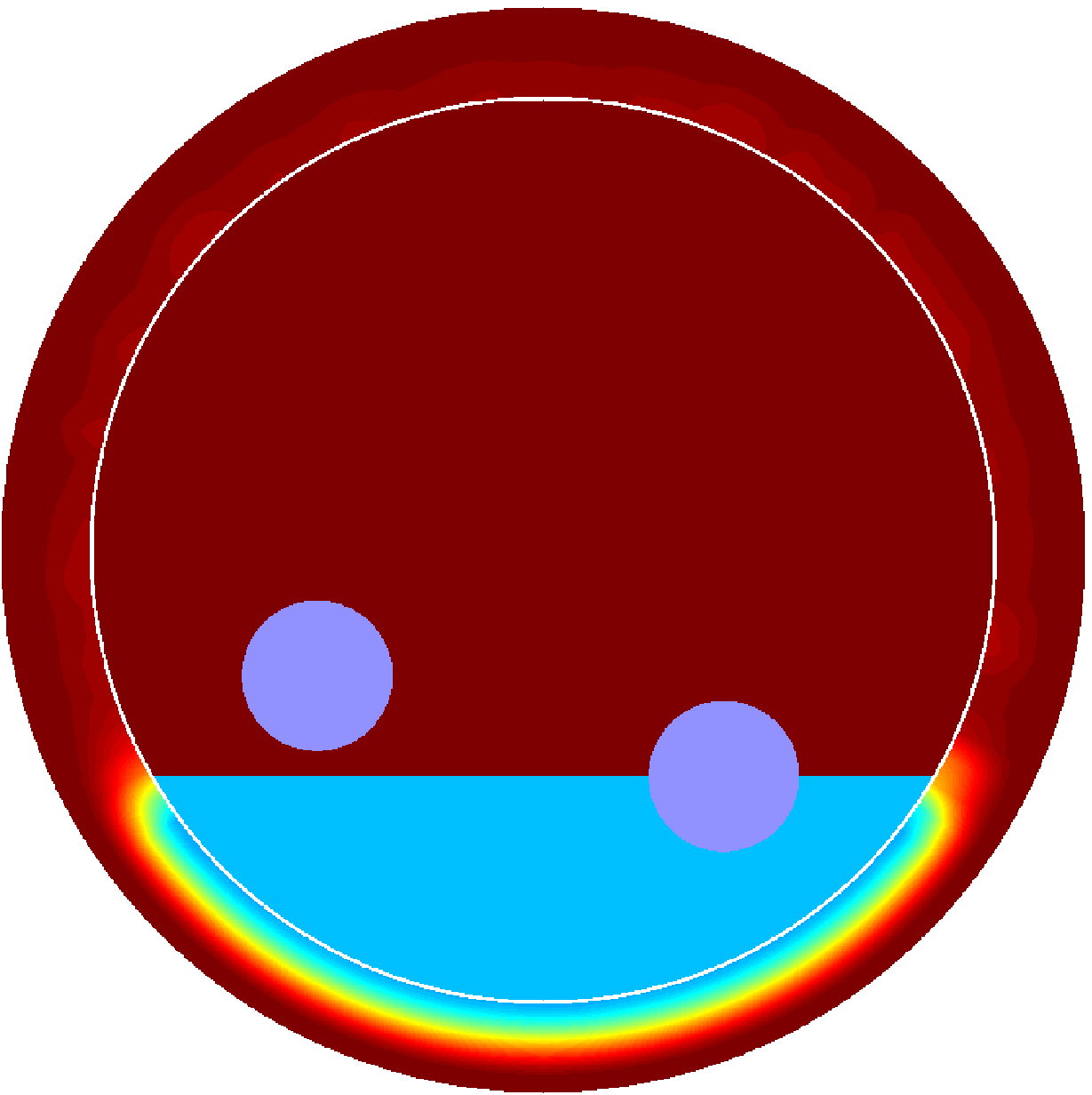}}
\epsfxsize=4.5cm
\put(0,47){\epsffile{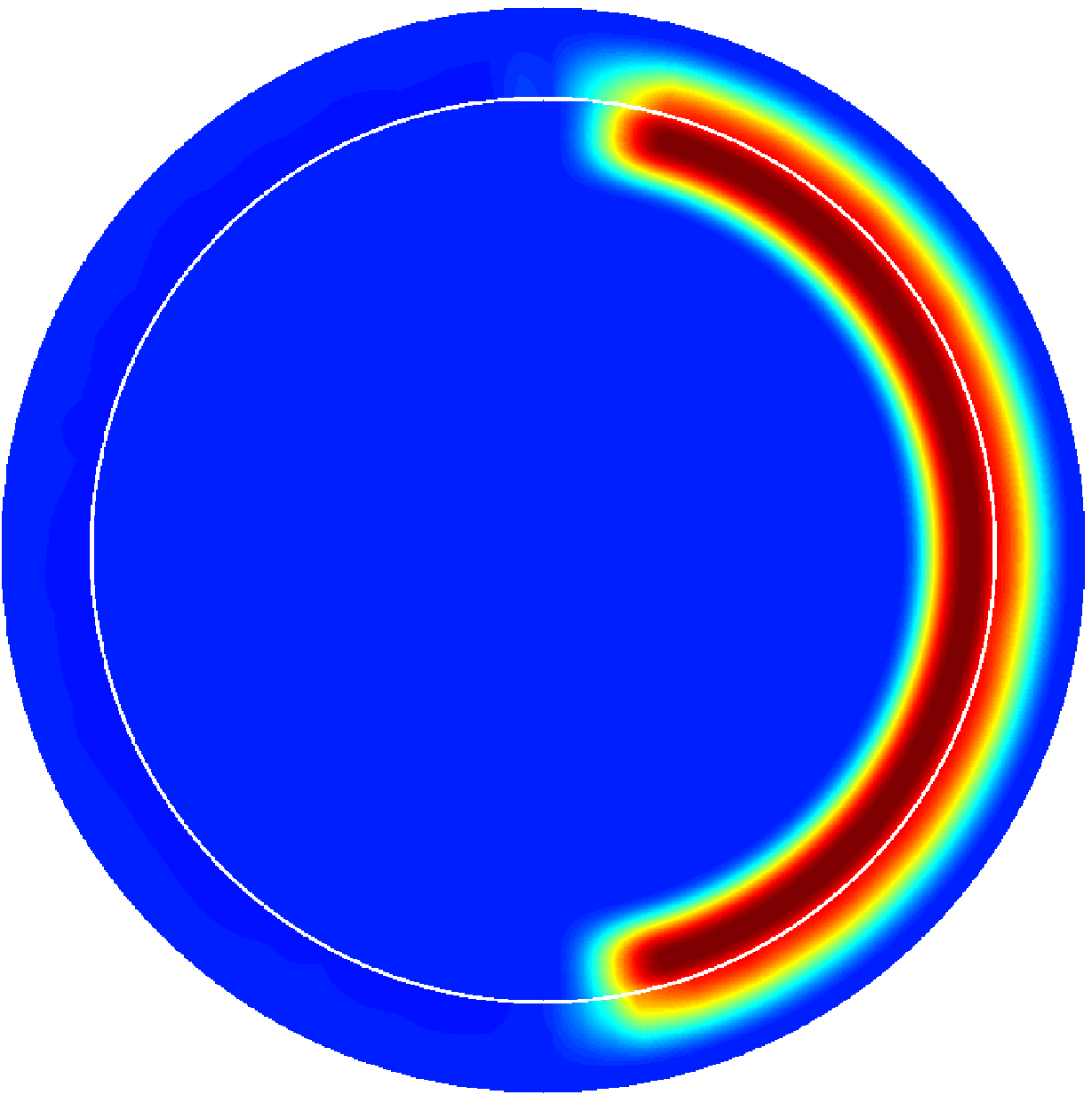}}
\epsfxsize=4.5cm
\put(0,94){\epsffile{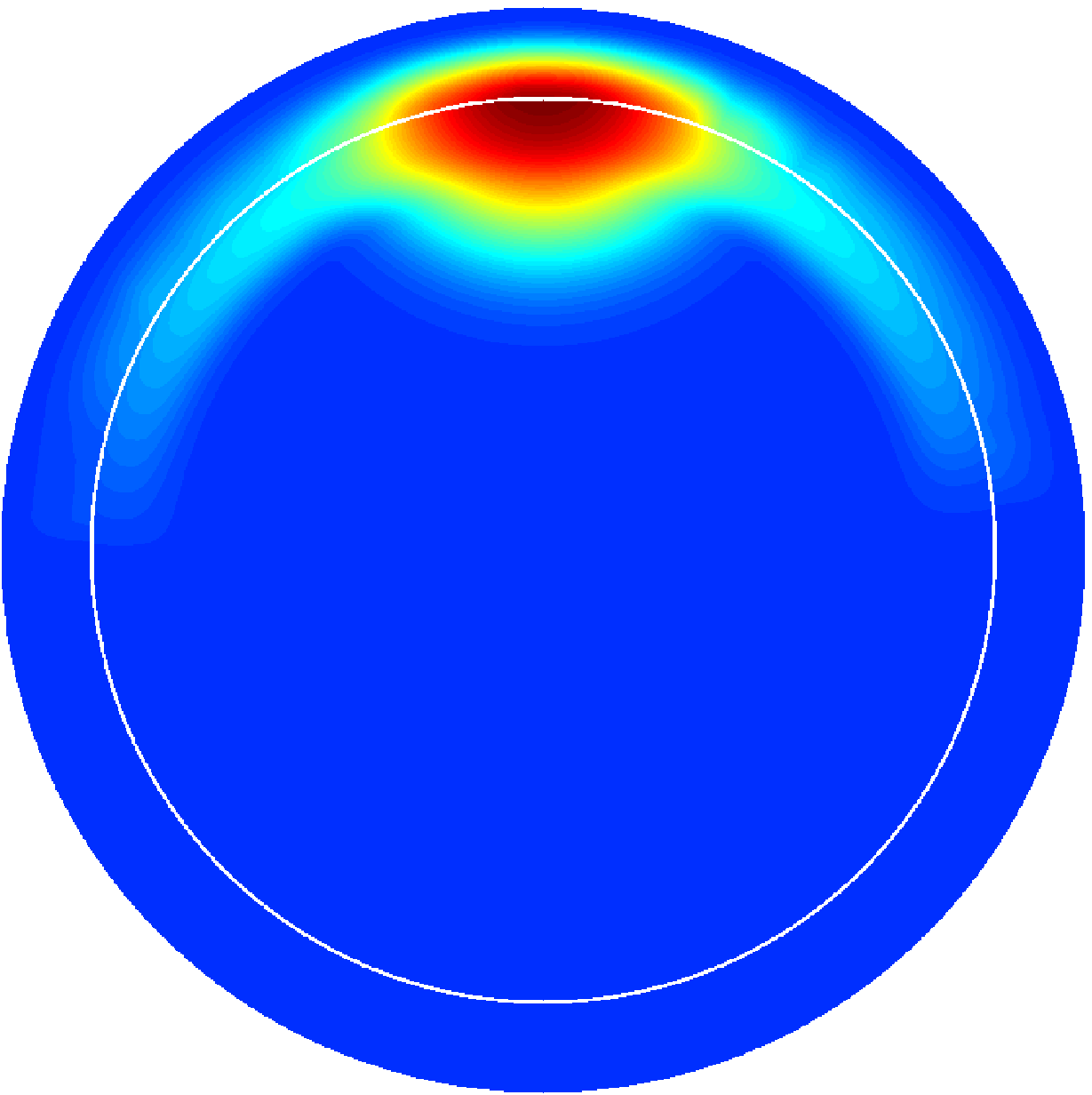}}
\epsfxsize=4.5cm
\put(0,141){\epsffile{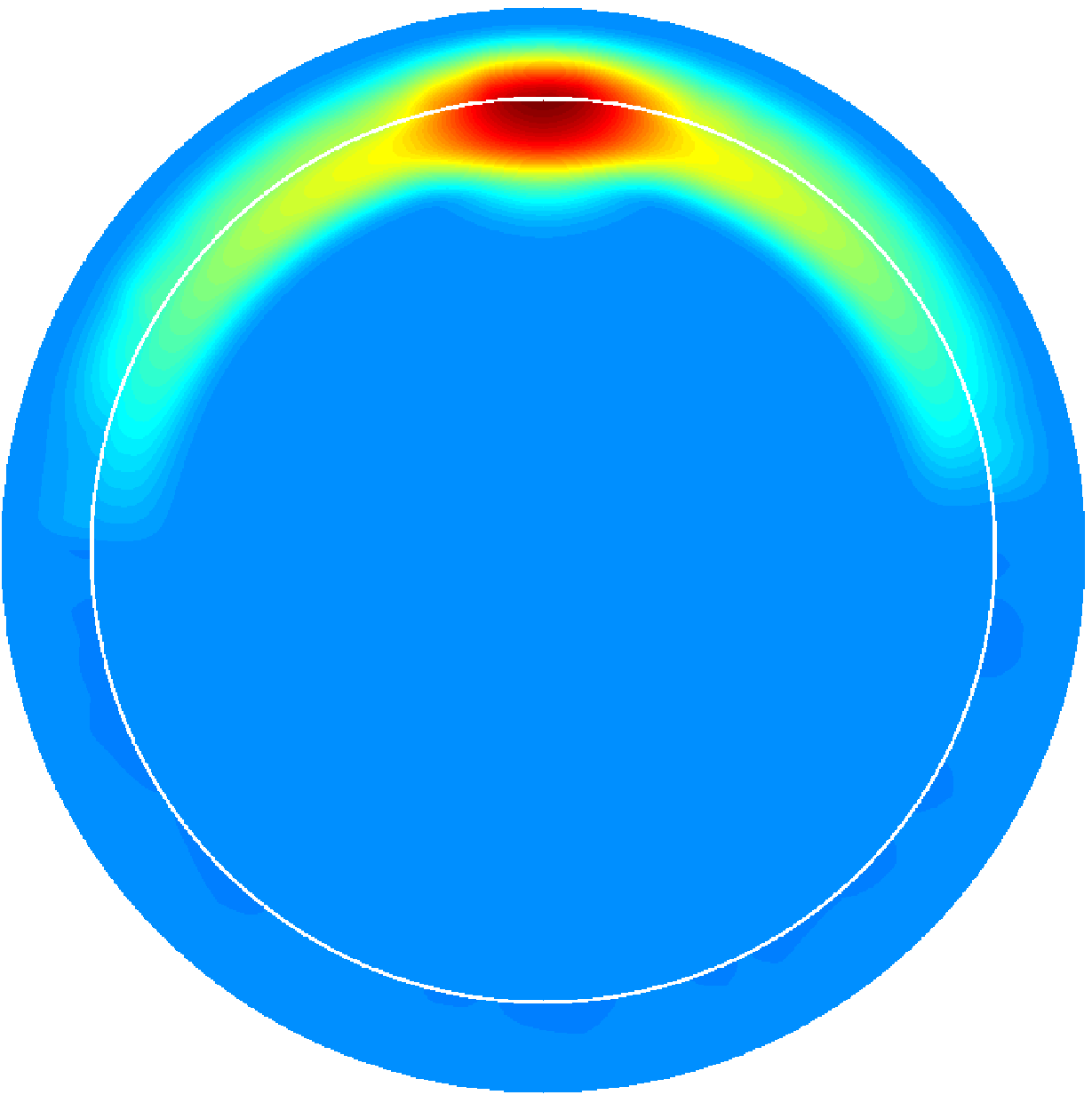}}
\epsfxsize=5cm
\put(65,2){\epsffile{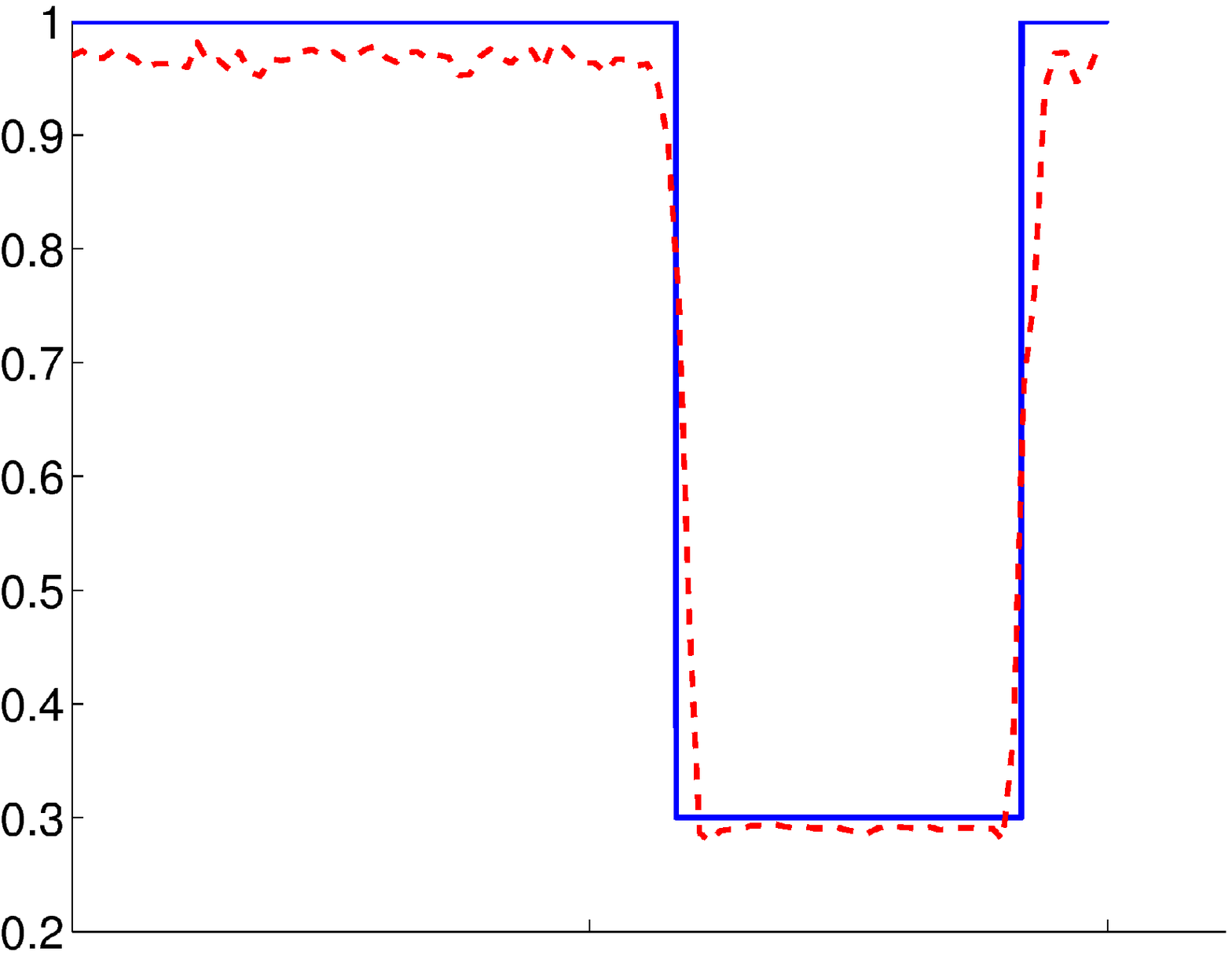}}
\epsfxsize=5cm
\put(65,50){\epsffile{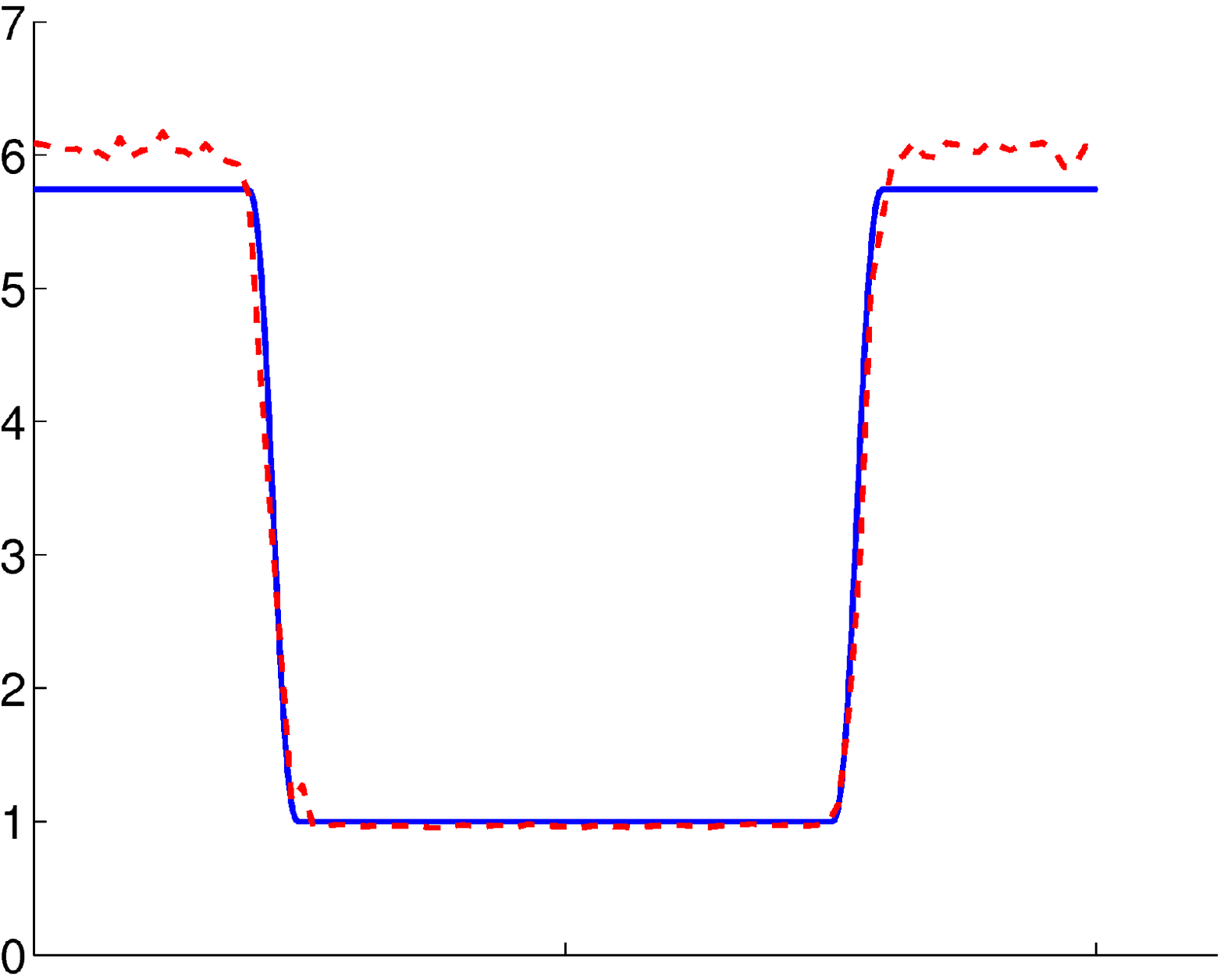}}
\epsfxsize=5cm
\put(65,98){\epsffile{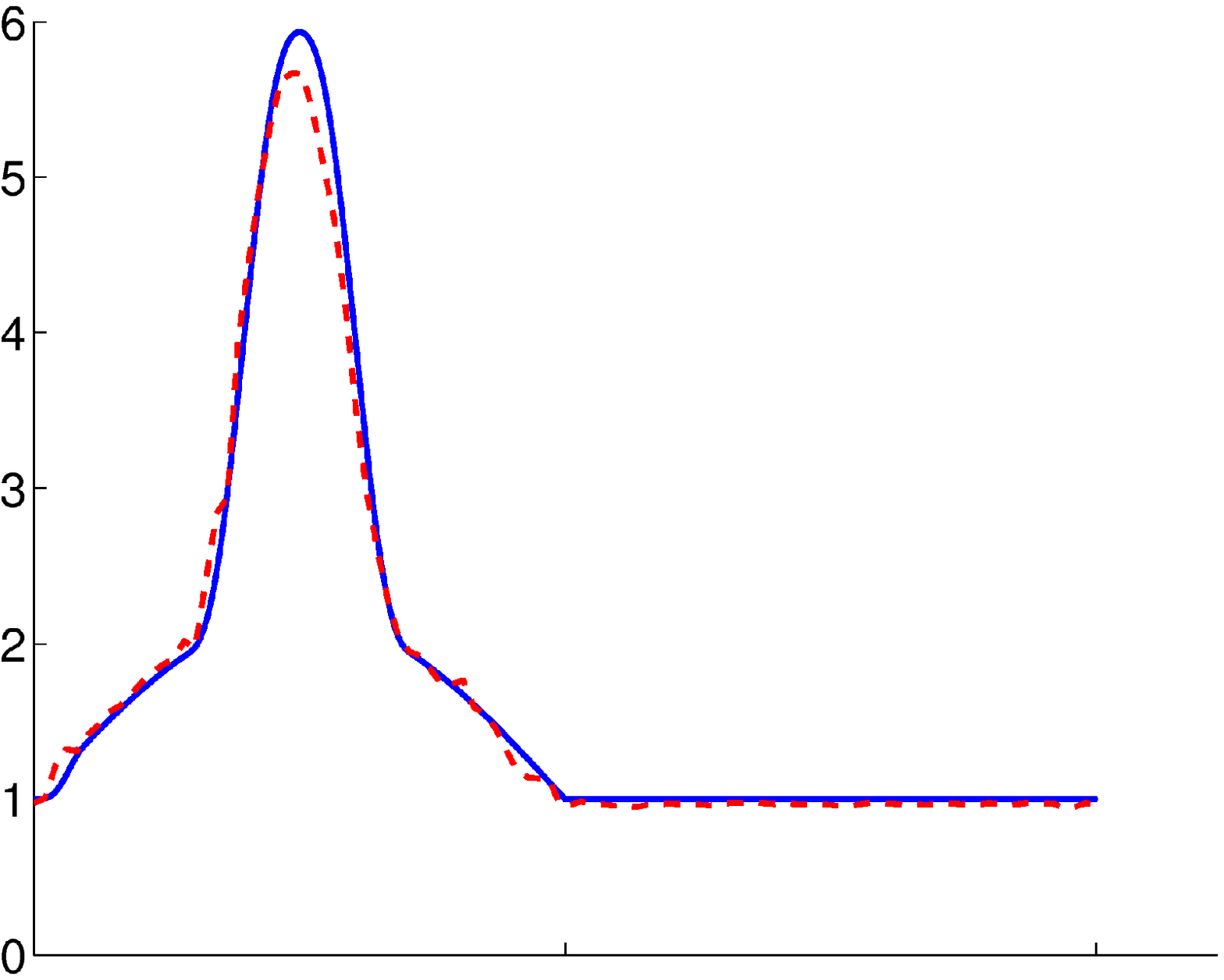}}
\epsfxsize=5cm
\put(65,146){\epsffile{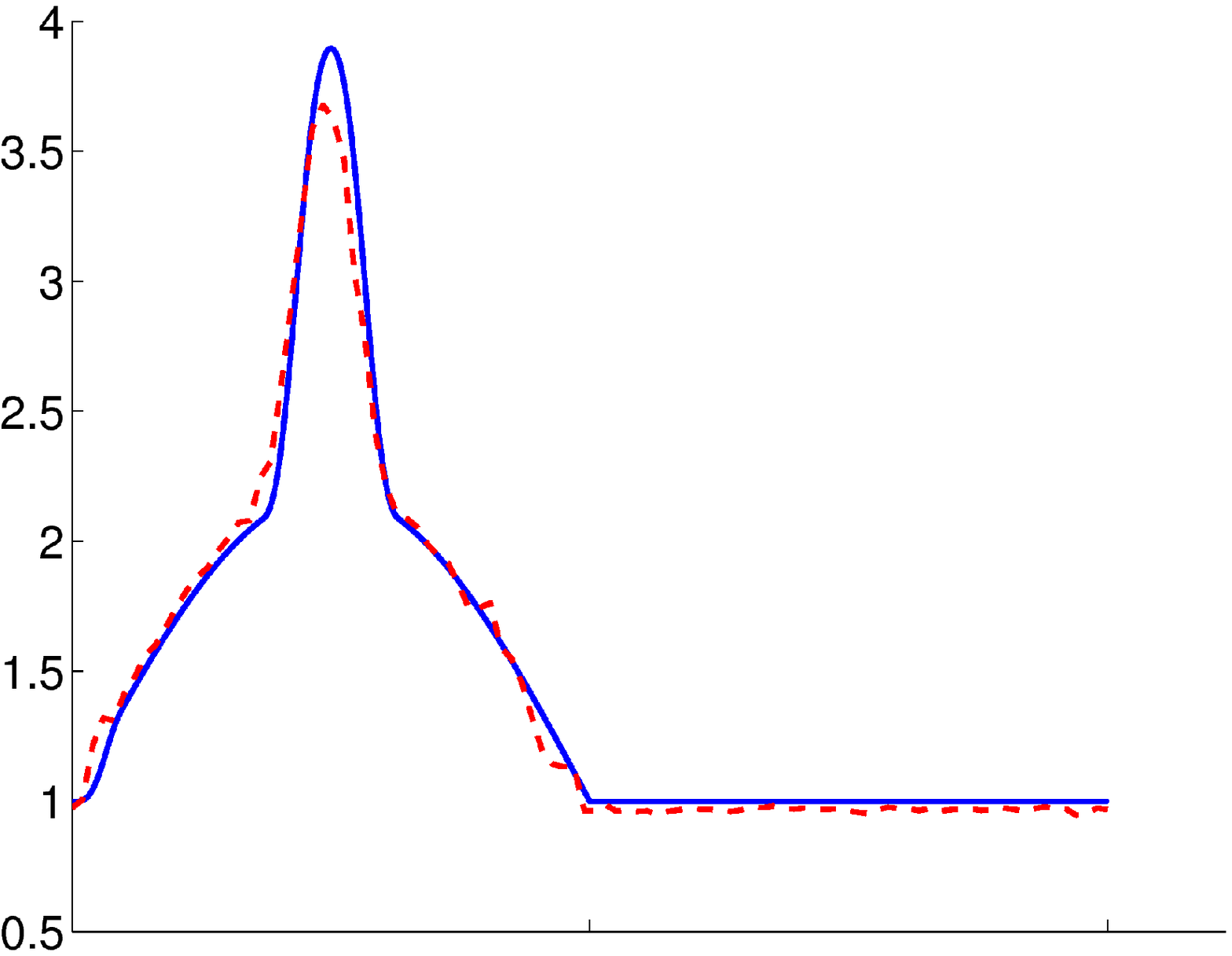}}
\put(40,183){ Example 1}
\put(40,136){ Example 2}
\put(40,89){ Example 3}
\put(40,42){ Example 4}
\put(97,143){\tiny $\theta$ over boundary}
\put(61,187){\tiny $\sigma|_{\bound_1}$}
\put(97,95){\tiny $\theta$ over boundary}
\put(61,139){\tiny $\sigma|_{\bound_1}$}
\put(97,47){\tiny $\theta$ over boundary}
\put(61,91){\tiny $\sigma|_{\bound_1}$}
\put(97,-1){\tiny $\theta$ over boundary}
\put(61,43){\tiny $\sigma|_{\bound_1}$}
\put(67,1){\tiny $0$}
\put(88,1){\tiny $\pi$}
\put(110,1){\tiny $2\pi$}
\put(67,49){\tiny $0$}
\put(88,49){\tiny $\pi$}
\put(110,49){\tiny $2\pi$}
\put(67,97){\tiny $0$}
\put(88,97){\tiny $\pi$}
\put(110,97){\tiny $2\pi$}
\put(67,145){\tiny $0$}
\put(88,145){\tiny $\pi$}
\put(110,145){\tiny $2\pi$}
\end{picture}
\caption{\label{fig:conductivities}Left column: Example conductivities
  $\gamma$ shown in the extended domain $\Om_2$, the white circle indicates the inner boundary $\bound_1$. Right column: actual
  traces of the conductivities at the inner boundary $\bound_1$ (solid
  line), and approximate traces at $\bound_1$ (dashed line) whose
  reconstruction is explained in Section \ref{sec:boundary}.}
\end{figure}

Based on the above numerical investigation we choose $r_2=1.2$ for the rest of this paper. We work with the following four example conductivities:
\begin{itemize}
\item Example one: conductivity has a high contrast bump right on the boundary 
$\bound_1$ and a circular inclusion near the boundary. All deviations
from background conductivity $1$ satisfy $\sigma(x)>1$.
\item Example two: similar to Example one but with a larger inclusion having higher conductivity.
\item Example three: conductivity has high-contrast behaviour near
  $\bound_1$, but the maximum of the deviation from background is not
  right at the boundary.  
\item Example four: crude model of a cross-section of an industrial pipeline, similar to the case in \cite{process_monitoring}. There is a sediment layer on the bottom of the tube, and two round low-conductivity inclusions.
\end{itemize}
See Figure \ref{fig:conductivities} for plots of the example conductivities and their traces on $\bound_1$.

We simulate noisy {\sc eit} data for each example conductivity using $c=0.00001$. The error $\norm{R^\epsilon_\sigma-R_\sigma}/\norm{R_\sigma}$ ranges between 0.00011 and 0.00076.

We use the method of Section \ref{sec:boundary} with
$M=32,\kappa=6,\alpha=4$ and 100 different angles to compute
approximately reconstructed traces $g$ on $\bound_1$. See the right column of Figure
\ref{fig:conductivities} for the result. Then, we compute the extended
conductivity $\gamma$ in the disc $\Om_2=D(0,r_2)=D(0,1.2)$ using
(\ref{sigma_tilde}) and (\ref{def:extension}) with the radius
$r_e=r_1+7/8(r_2-r_1) = 1.175$. Since $g$ is only approximately the
same as $\sigma|_{\bound_1}$ there are discontinuities in $\gamma$ in
all cases.


We compute the intermediate {\sc dn} maps $\Lambda^{ij}$ using the finite
element method and 425984 triangles in the annulus $\diff$.
To check the accuracy of formula (\ref{L_gamma}) we also calculate
$\Lambda_\gamma$ directly by the finite element method (and 1081344 triangles in $\Om_2$) using the knowledge of
$\gamma$. The error $\norm{L^\epsilon_\gamma-L^2_\gamma}/\norm{L^2_\gamma}$,
where $L^2_\gamma$ is the {\sc dn} map calculated directly on the boundary
$\bound_2$, was less than 2.2\% in all cases. The condition number of the matrix $L_\sigma^\ep-L^{11}$ used in \eqref{L_gamma} was less than $27$ in all test cases.

\begin{figure}[p]
\begin{picture}(120,110)
\epsfxsize=13cm
\put(0,0){\epsffile{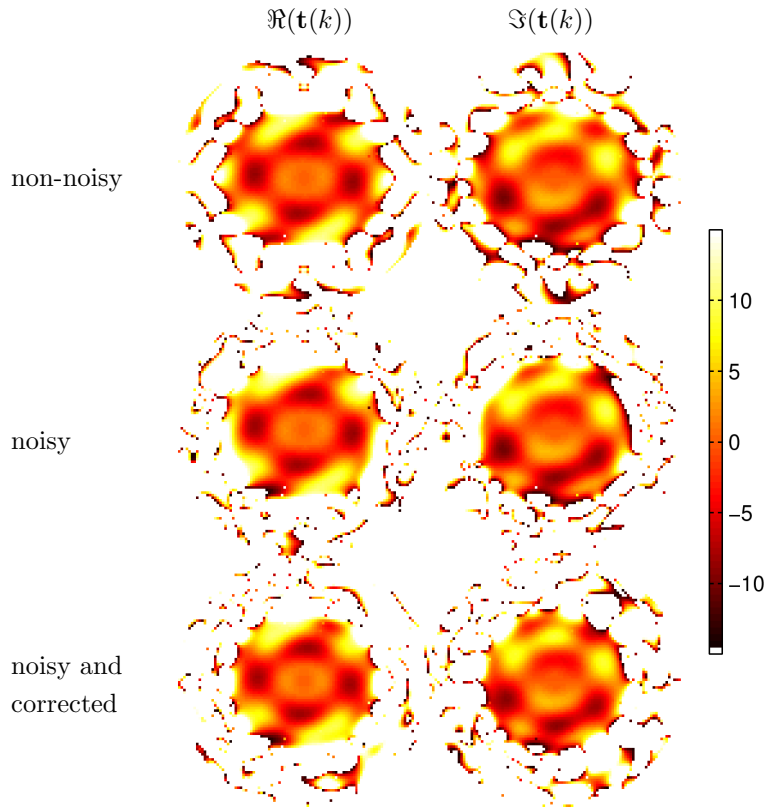}}
\put(43,106){$\Re(\T(k))$}
\put(75,106){$\Im(\T(k))$}
\put(9,85){non-noisy}
\put(9,50){noisy}
\put(9,20){noisy and}
\put(9,15){corrected}
\end{picture}
\caption{\label{fig:scattering_transforms} The scattering transform in example four. The first row shows the the real and imaginary parts of \eqref{eq:Teps} substituting $L_\sigma$ in place of $\L_\sigma^\ep$. The second row shows the same functions using $L_\sigma^\ep$, and the third row is again the same, but uses $L_\gamma^\ep$ calculated from \eqref{L_gamma}. The real part of $\T_R(k)$ is in the left column, the imaginary part on the right. The scattering transform is calculated in a grid of spectral parameters $k$, where $\abs{k}<10$. In white areas we have $\abs{\T_R(k)}>15$, meaning the calculation has failed or is close to failing due to computational error caused by large values of $\abs{k}$. }
\end{figure}

Figure \ref{fig:scattering_transforms} illustrates how the noise and the boundary correction procedure affect the scattering transform in example four. The first row shows the the real and imaginary parts of \eqref{eq:Teps} substituting $L_\sigma$ in place of $L_\sigma^\ep$. The second row shows the same functions using $L_\sigma^\ep$, and the third row is again the same, but uses $L_\gamma^\ep$ calculated from \eqref{L_gamma}. The real part of $\T_R(k)$ is in the left column, the imaginary part on the right. The scattering transform is calculated in a grid of spectral parameters $k$, where $\abs{k}<10$. In white areas we have $\abs{\T_R(k)}>15$, meaning the calculation has failed or is close to failing due to computational error caused by large values of $\abs{k}$.

\begin{figure}[p]
\begin{picture}(110,92)
\epsfxsize=5.5cm
\put(0,47){\epsffile{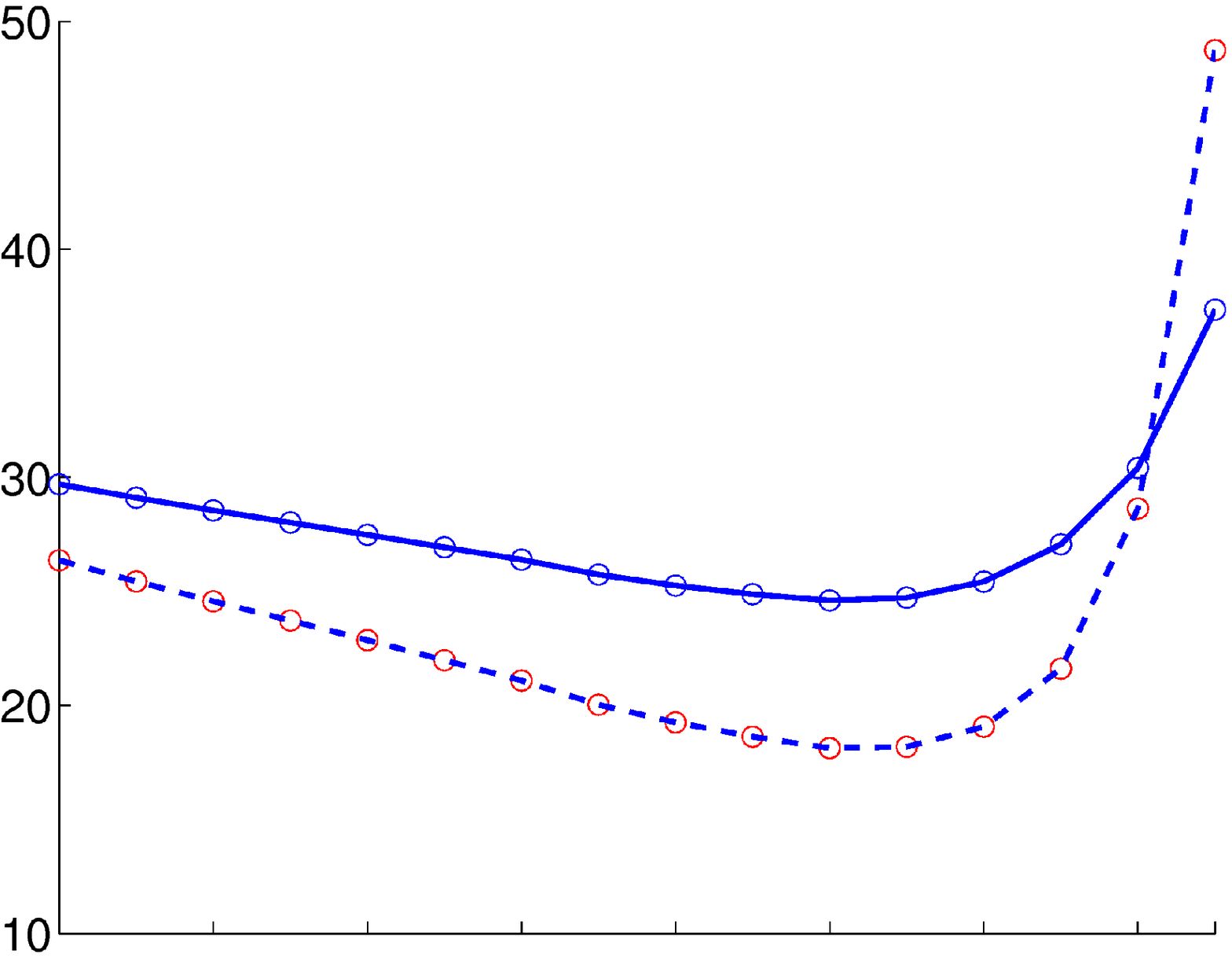}}
\epsfxsize=5.5cm
\put(59,47){\epsffile{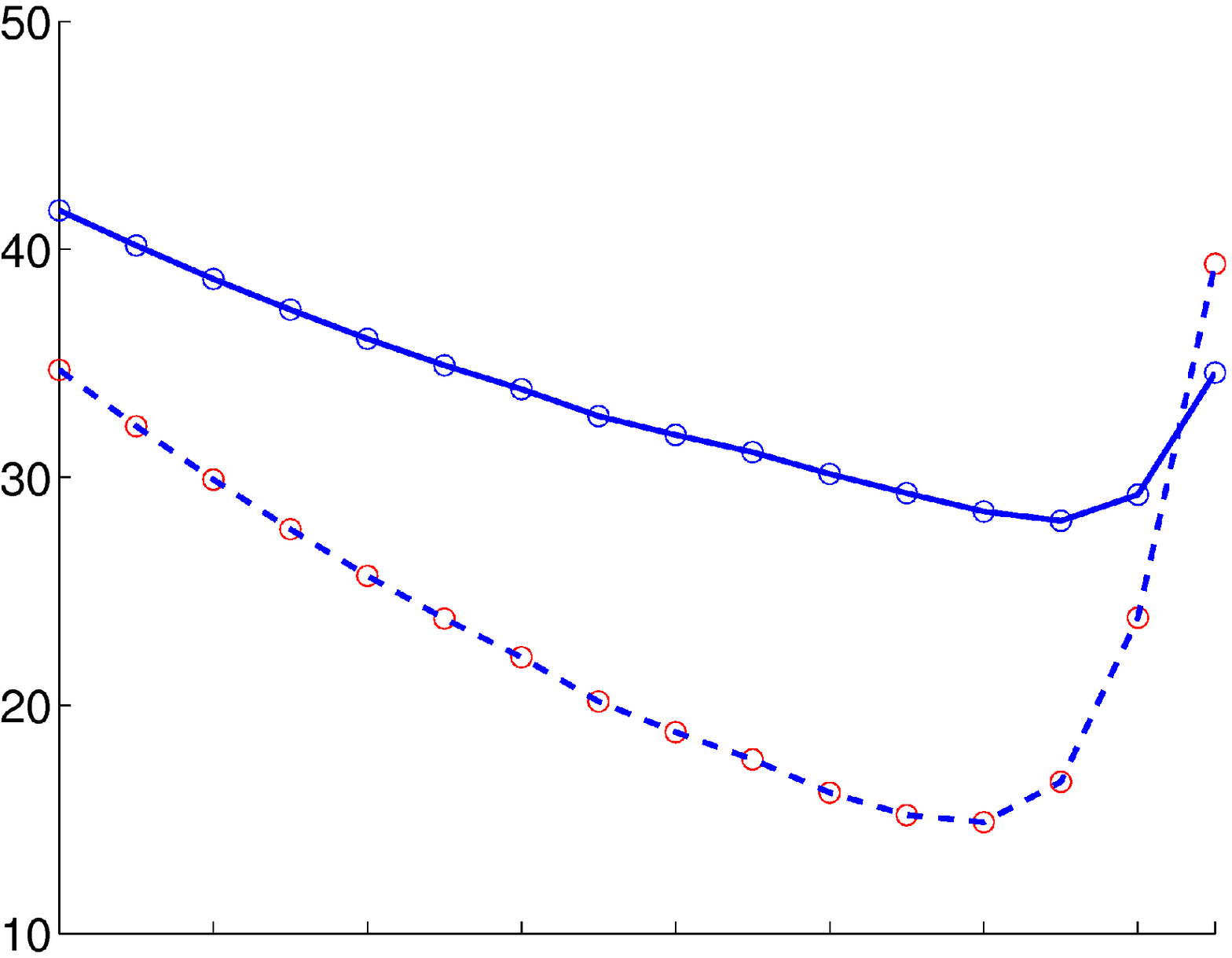}}
\epsfxsize=5.5cm
\put(0,0){\epsffile{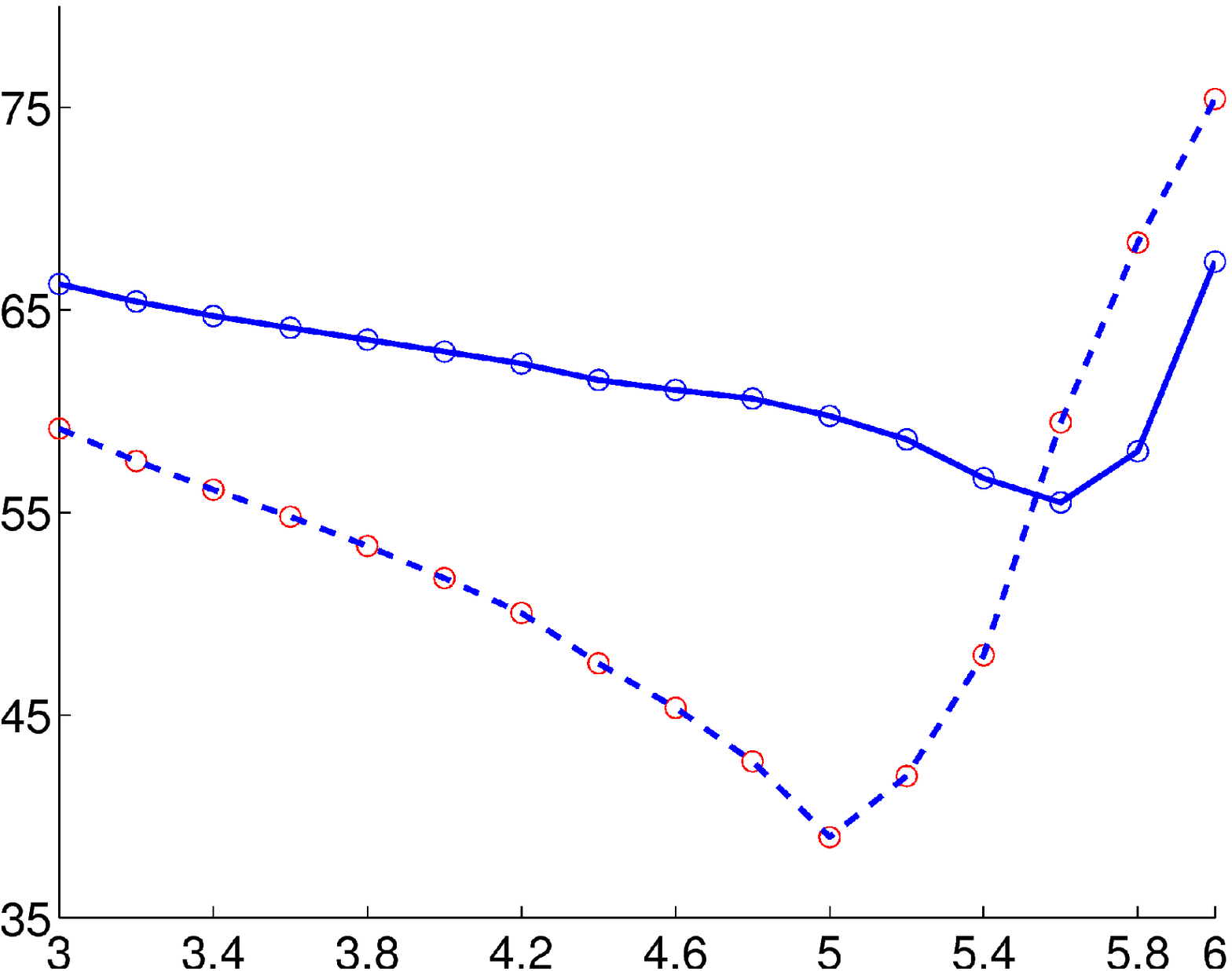}}
\epsfxsize=5.5cm
\put(59,0){\epsffile{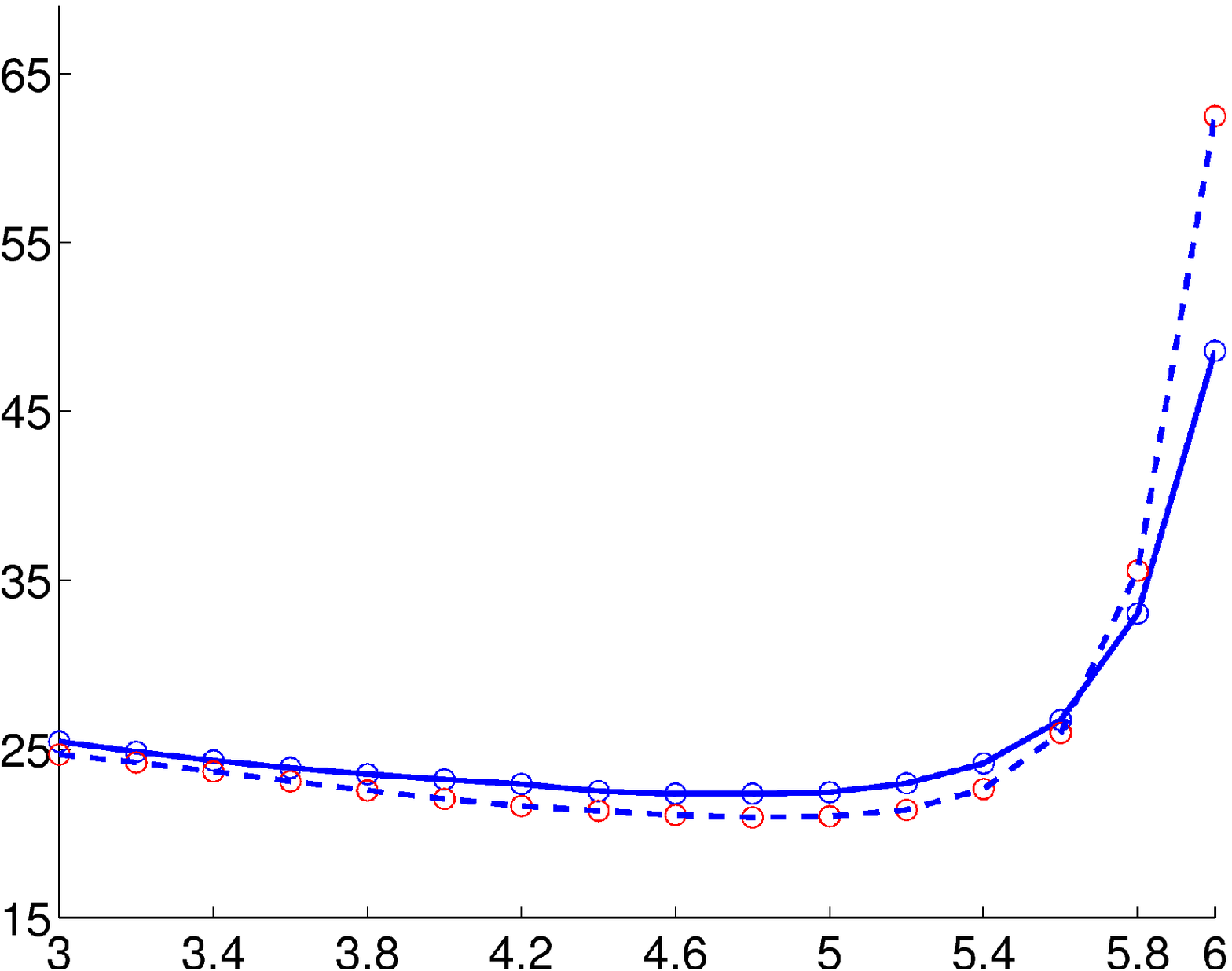}}
\put(-2,85){\tiny \%}
\put(52,-2){\tiny $R$}
\put(-2,42){\tiny \%}
\put(111,-2){\tiny $R$}
\put(7,85){ \small Example 1}
\put(67,85){\small Example 2}
\put(7,39){\small Example 3}
\put(67,39){\small Example 4}
\end{picture}
\caption{\label{fig:errors}$L^2$-error graphs as a function of truncation radius $R$ of the scattering transform for different examples;
  solid line is for the traditional D-bar reconstructions, dashed line is for
  boundary corrected reconstructions. The $R$ -axis is the same in all four plots.}
\end{figure}

\begin{figure}[p]
\begin{picture}(120,86)
\epsfxsize=3cm
\put(0,24){\epsffile{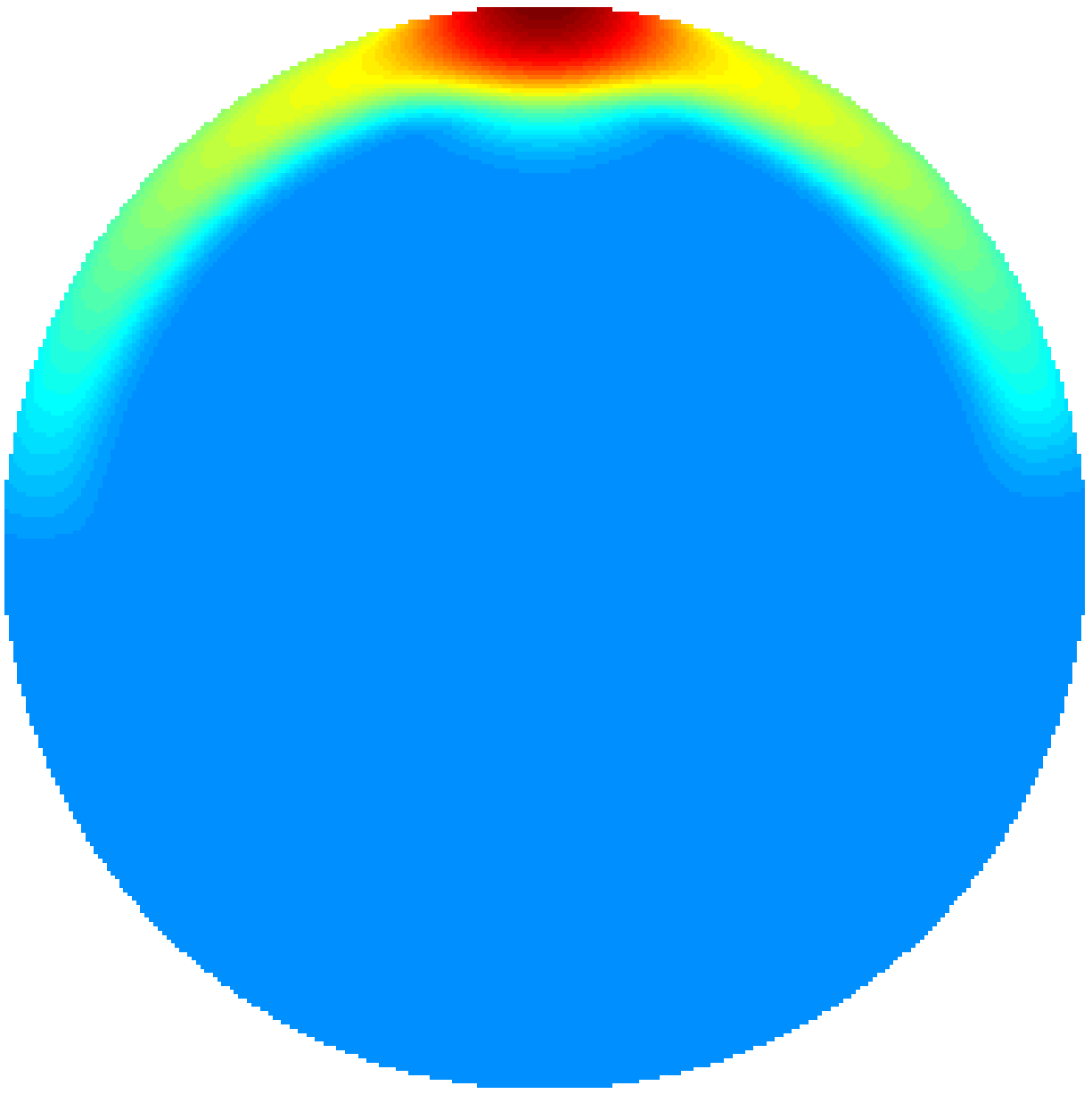}}
\epsfxsize=3cm
\put(30,48){\epsffile{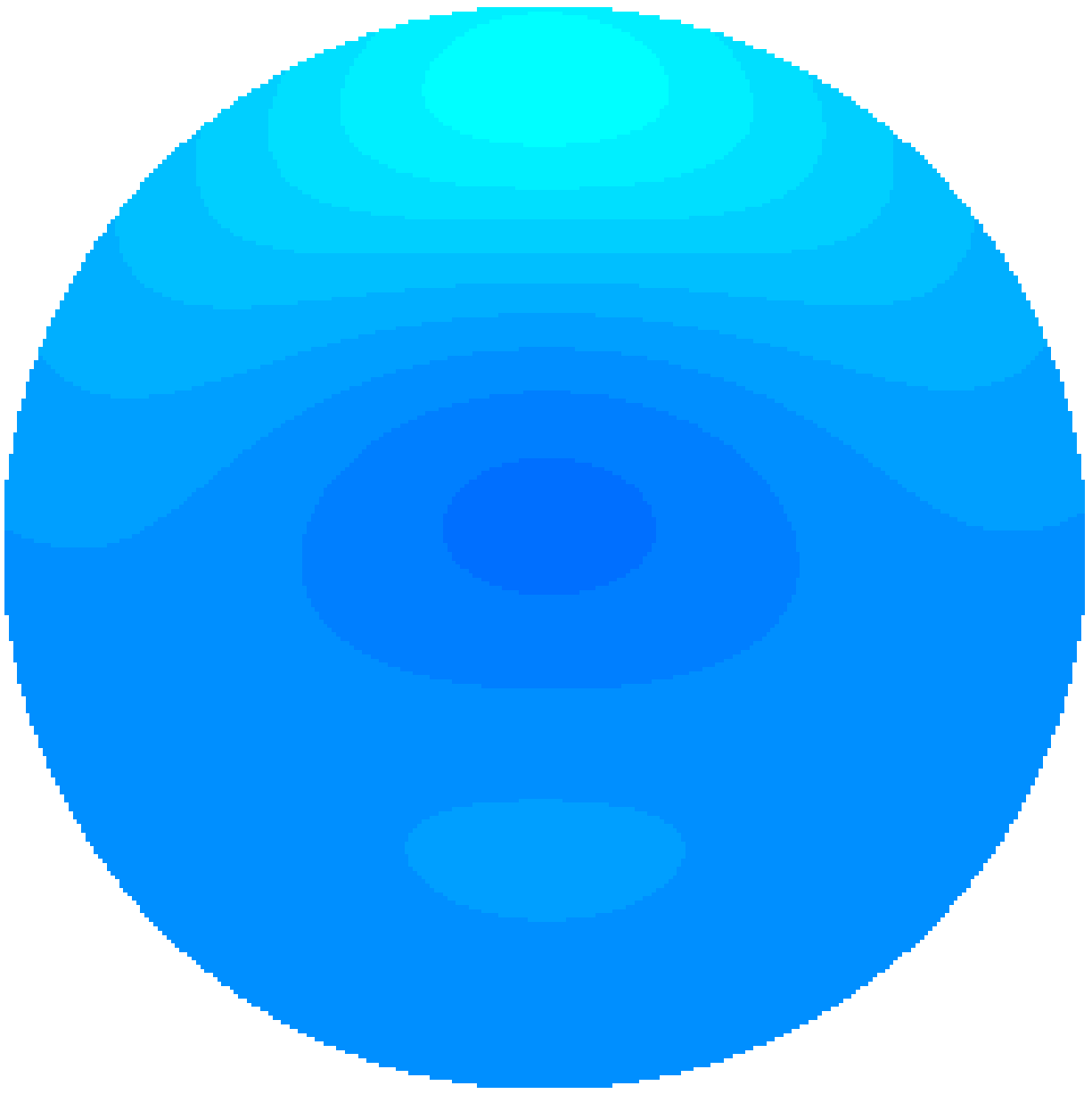}}
\epsfxsize=3cm
\put(60,48){\epsffile{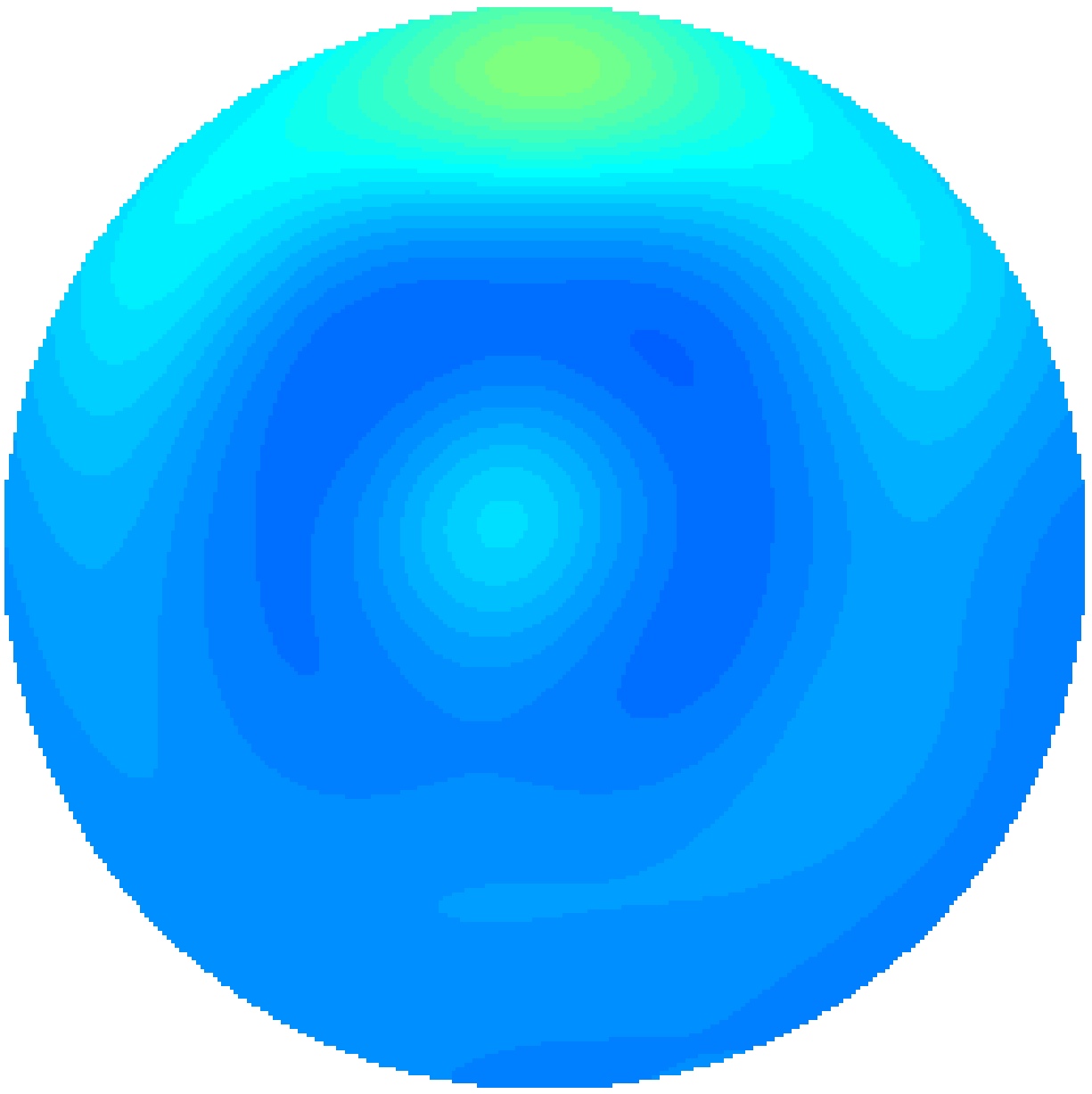}}
\epsfxsize=3cm
\put(90,48){\epsffile{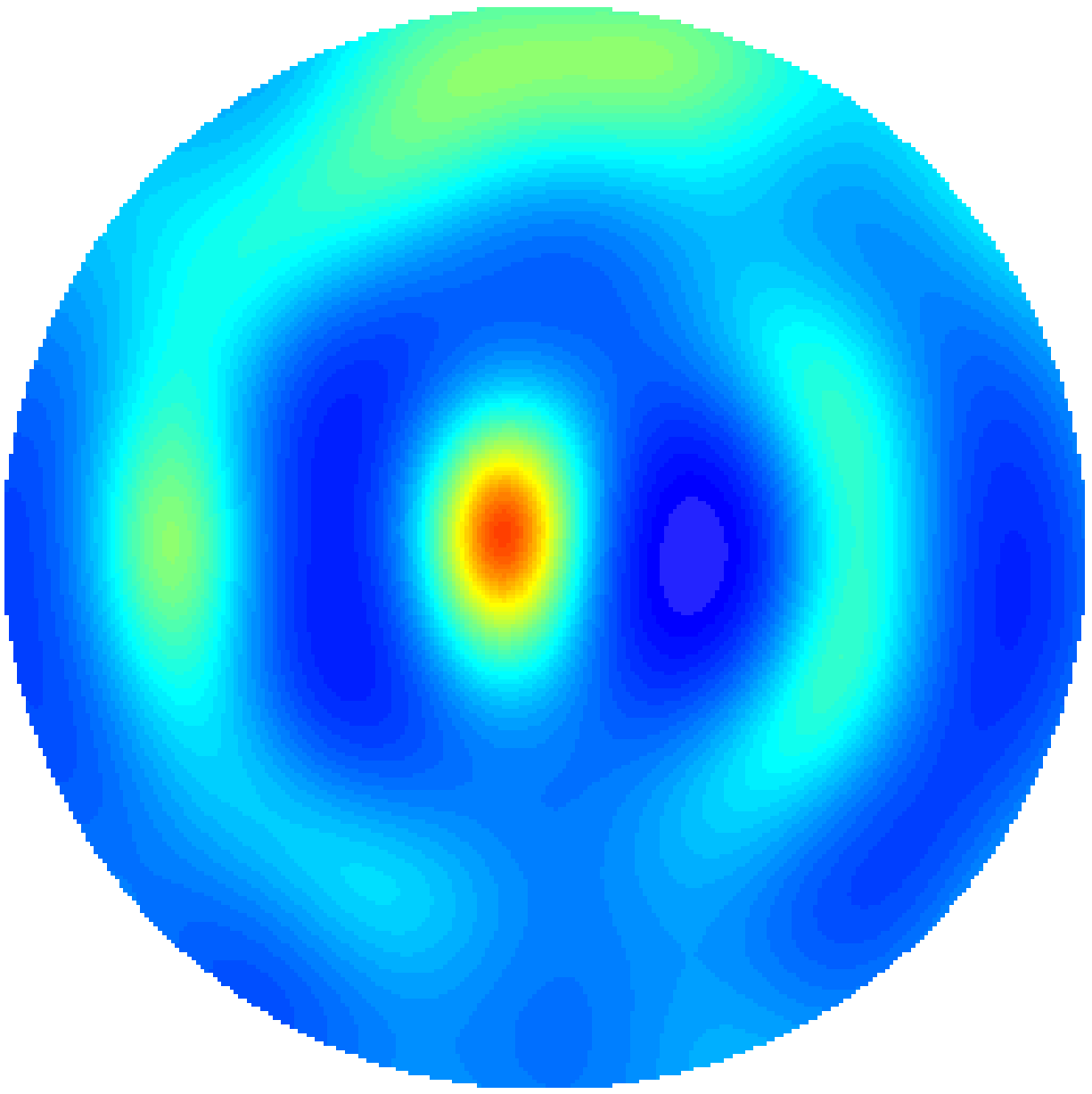}}
\epsfxsize=3cm
\put(30,0){\epsffile{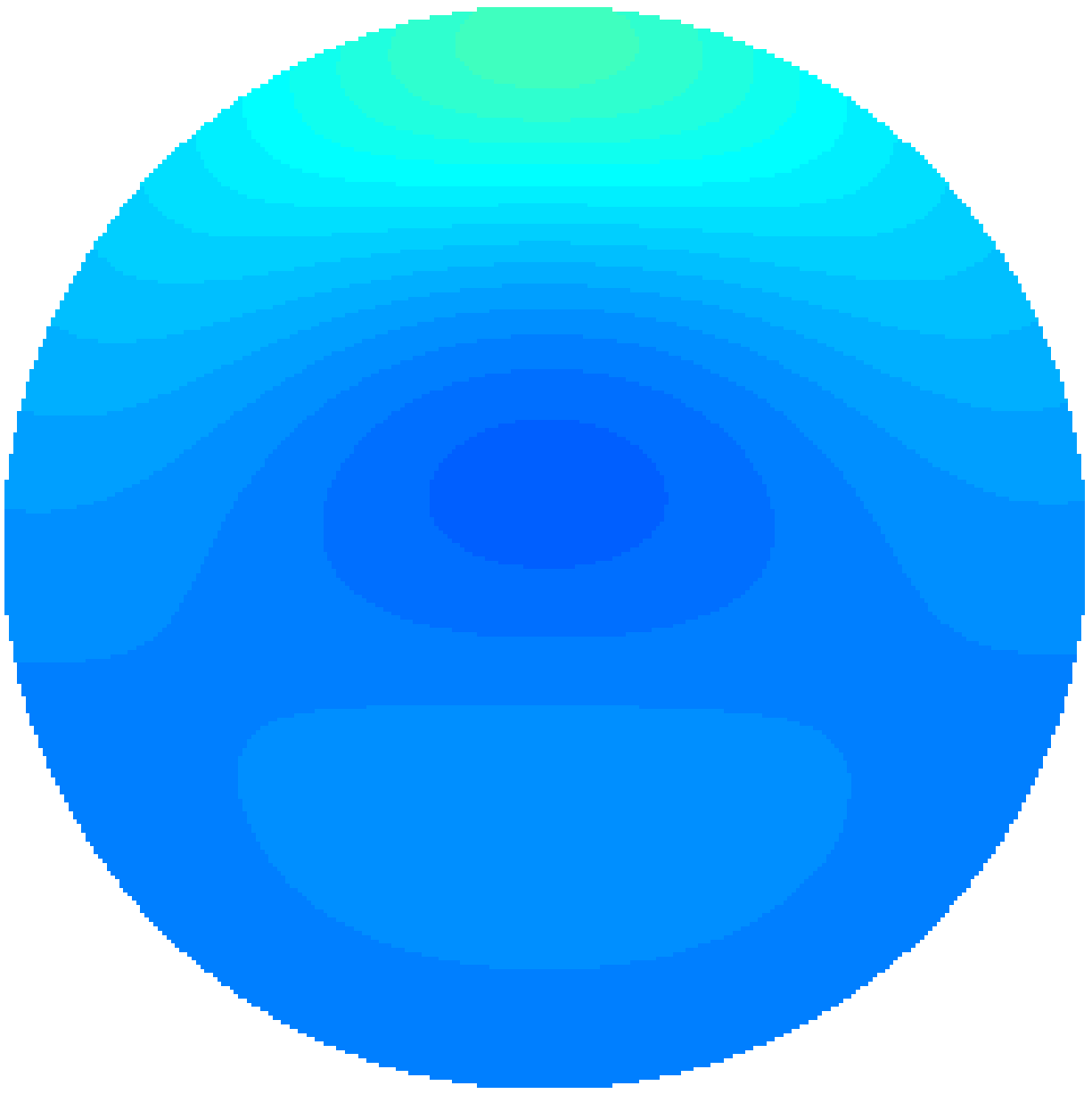}}
\epsfxsize=3cm
\put(60,0){\epsffile{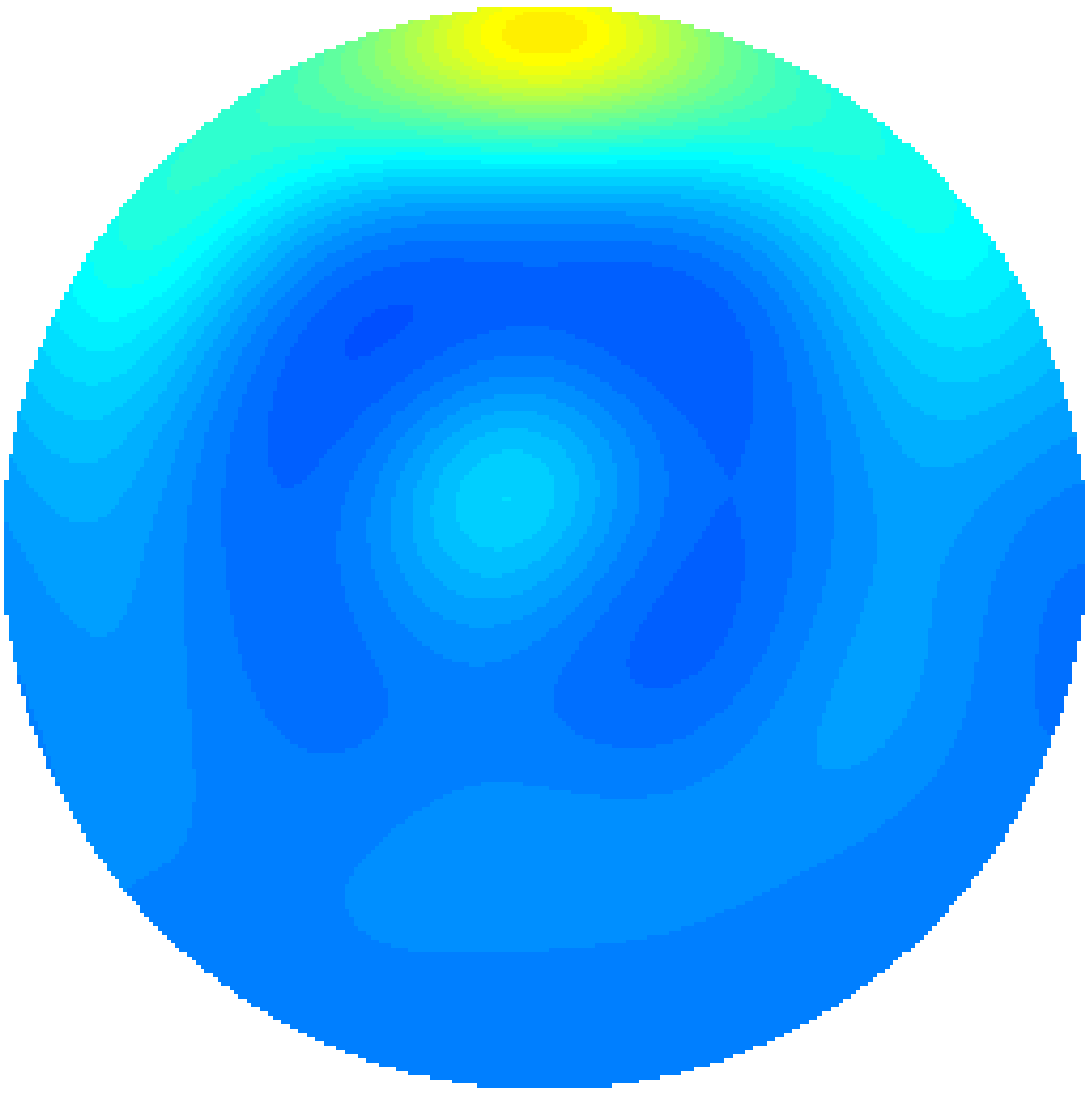}}
\epsfxsize=3cm
\put(90,0){\epsffile{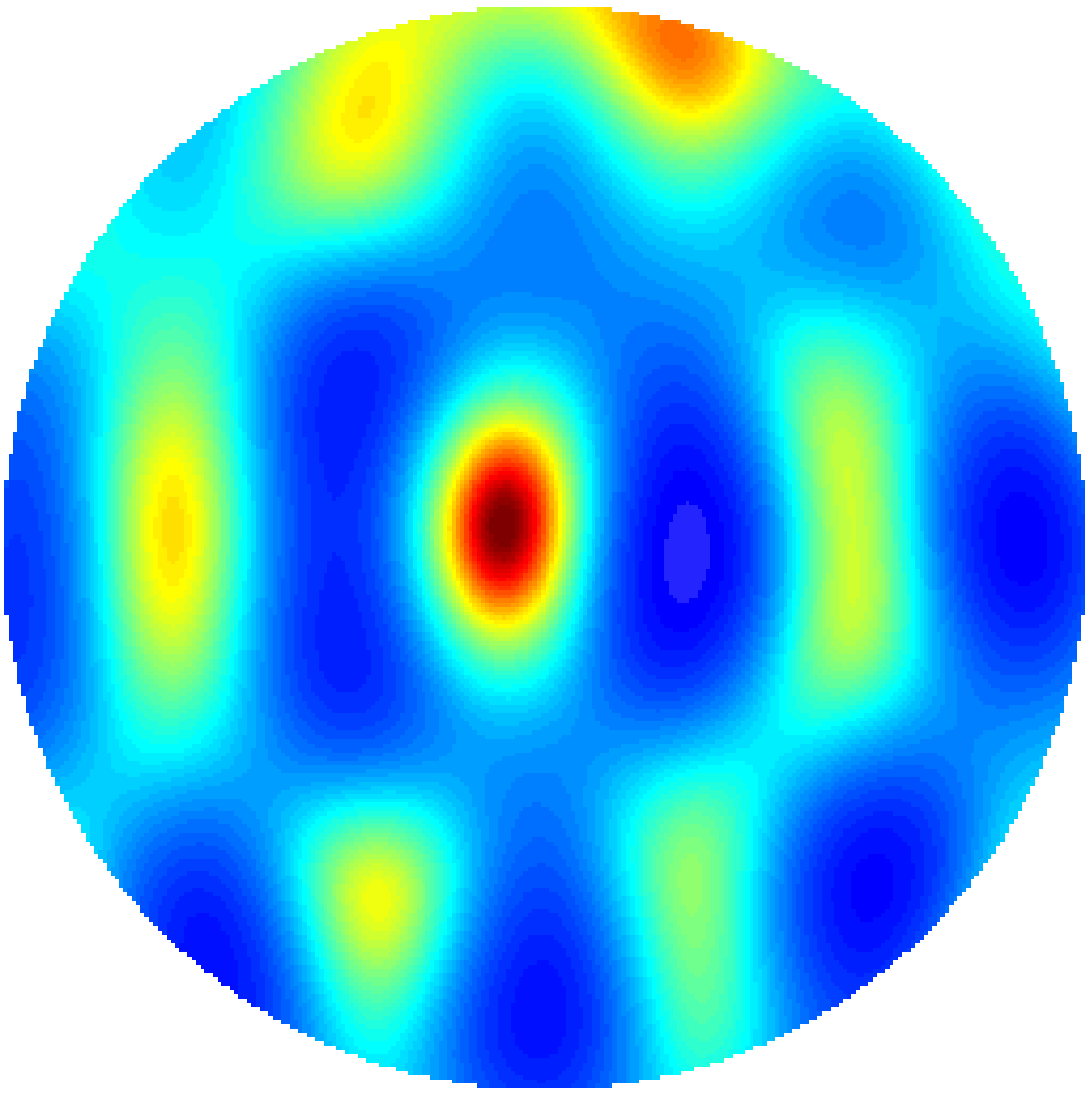}}

\put(6,59){\small Example 1}
\put(-1,55){\small Original conductivity}

\put(37,83){\small $R=3.0$}
\put(35,79){\small Uncorrected}
\put(67,83){\small $R=5.0$}
\put(65,79){\small Uncorrected}
\put(97,83){\small $R=6.0$}
\put(95,79){\small Uncorrected}

\put(37,35){\small $R=3.0$}
\put(37,31){\small Corrected}
\put(67,35){\small $R=5.0$}
\put(67,31){\small Corrected}
\put(97,35){\small $R=6.0$}
\put(97,31){\small Corrected}

\put(30,47){30\%}
\put(60,47){25\%}
\put(90,47){37\%}
\put(30,-1){26\%}
\put(60,-1){18\%}
\put(90,-1){49\%}

\end{picture}
\caption{\label{fig:ex1recon}Example 1 reconstructions; the original
  conductivity on the left, traditional D-bar reconstructions on the upper
  row and boundary corrected reconstructions on the lower row; the numbers
  beside the pictures are $L^2$ -errors, for the full error graph, see figure
  \ref{fig:errors}. The first reconstruction pair is always calculated with
$R=3$, the second one is the one with the lowest numerical $\L^2$
-error for the boundary corrected reconstruction, and the third one is with $R=6$ to show how the reconstructions fail.}
\end{figure}

\begin{figure}[p]
\begin{picture}(120,86)
\epsfxsize=3cm
\put(0,24){\epsffile{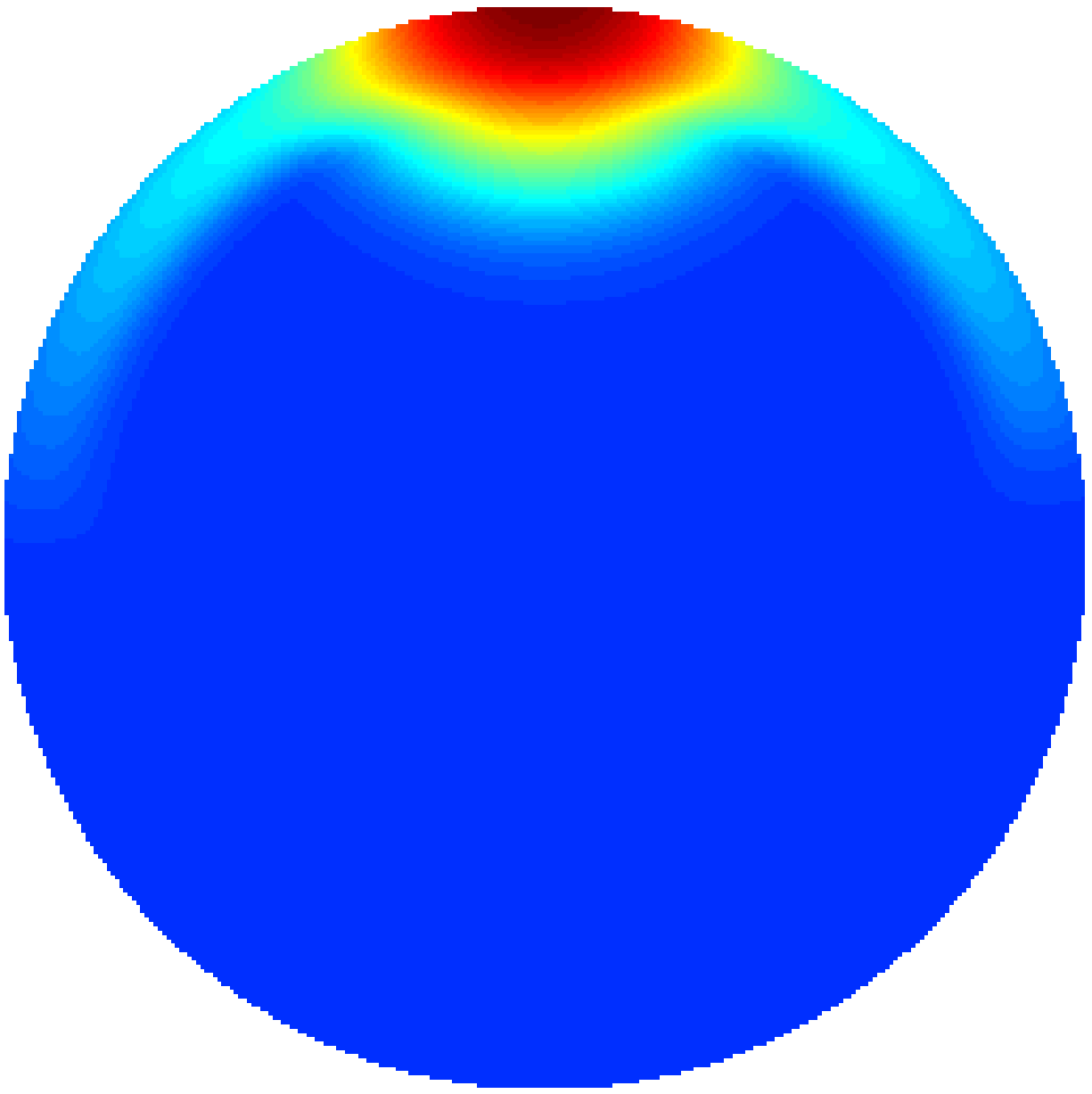}}
\epsfxsize=3cm
\put(30,48){\epsffile{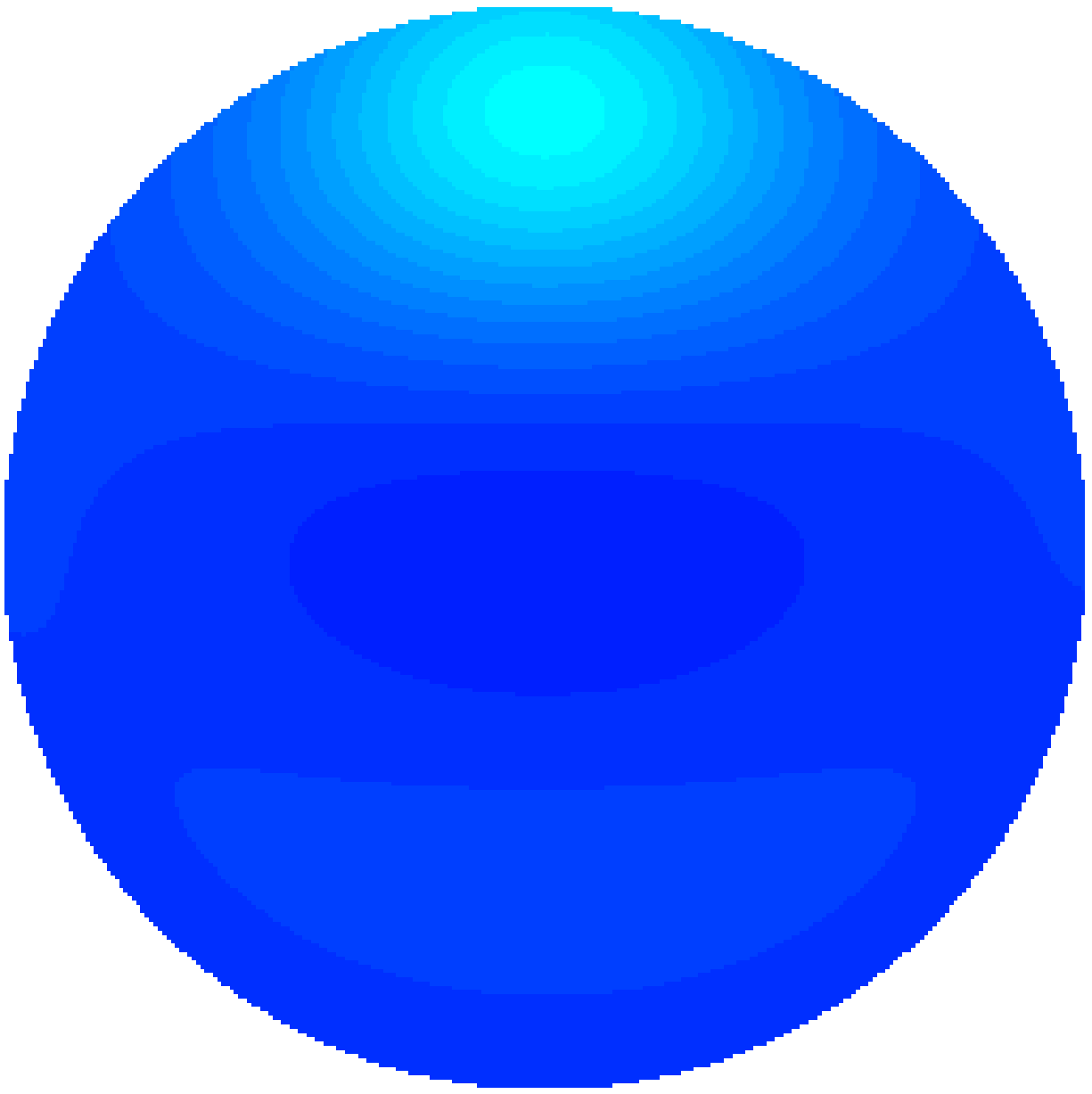}}
\epsfxsize=3cm
\put(60,48){\epsffile{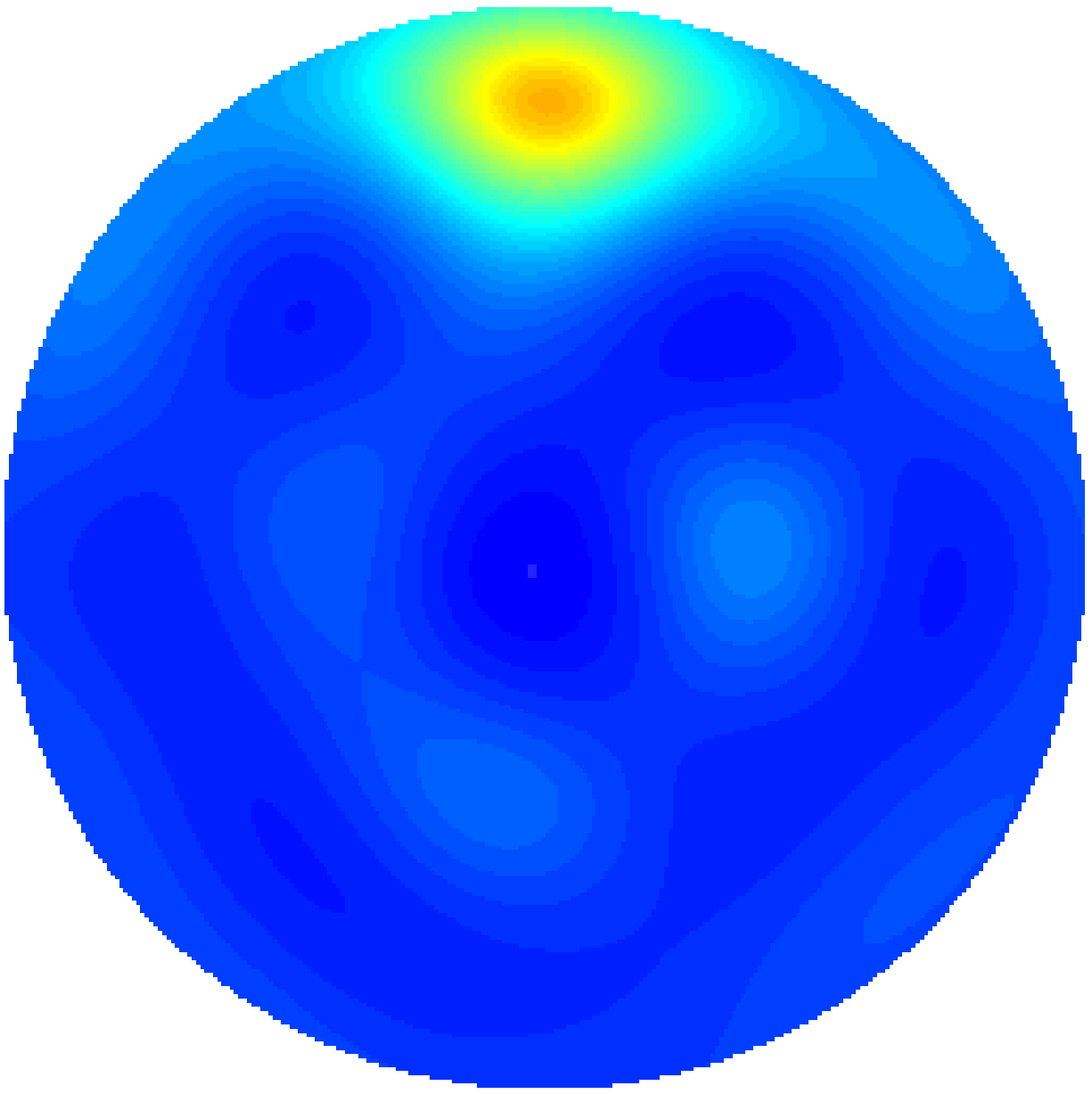}}
\epsfxsize=3cm
\put(90,48){\epsffile{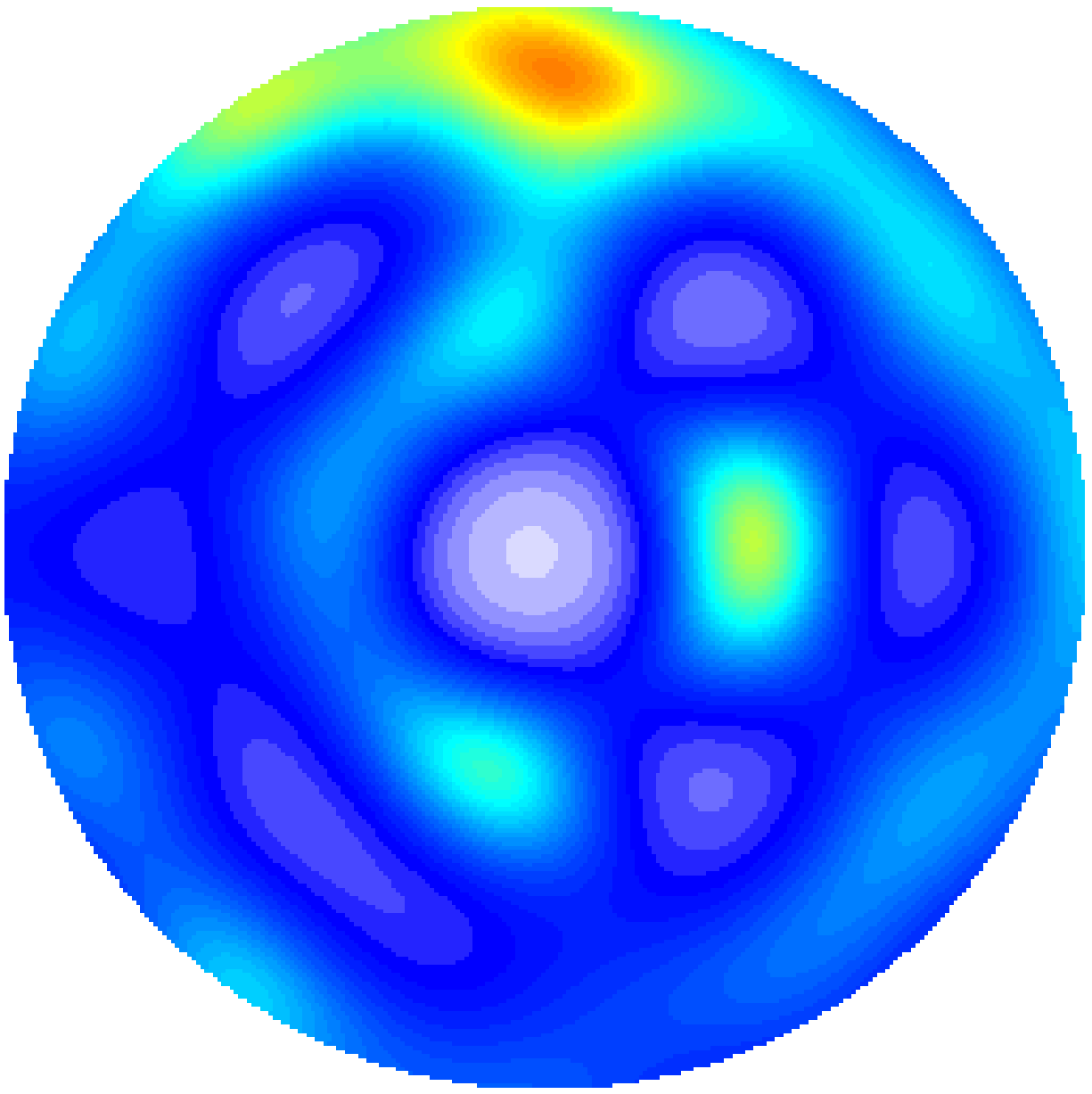}}
\epsfxsize=3cm
\put(30,0){\epsffile{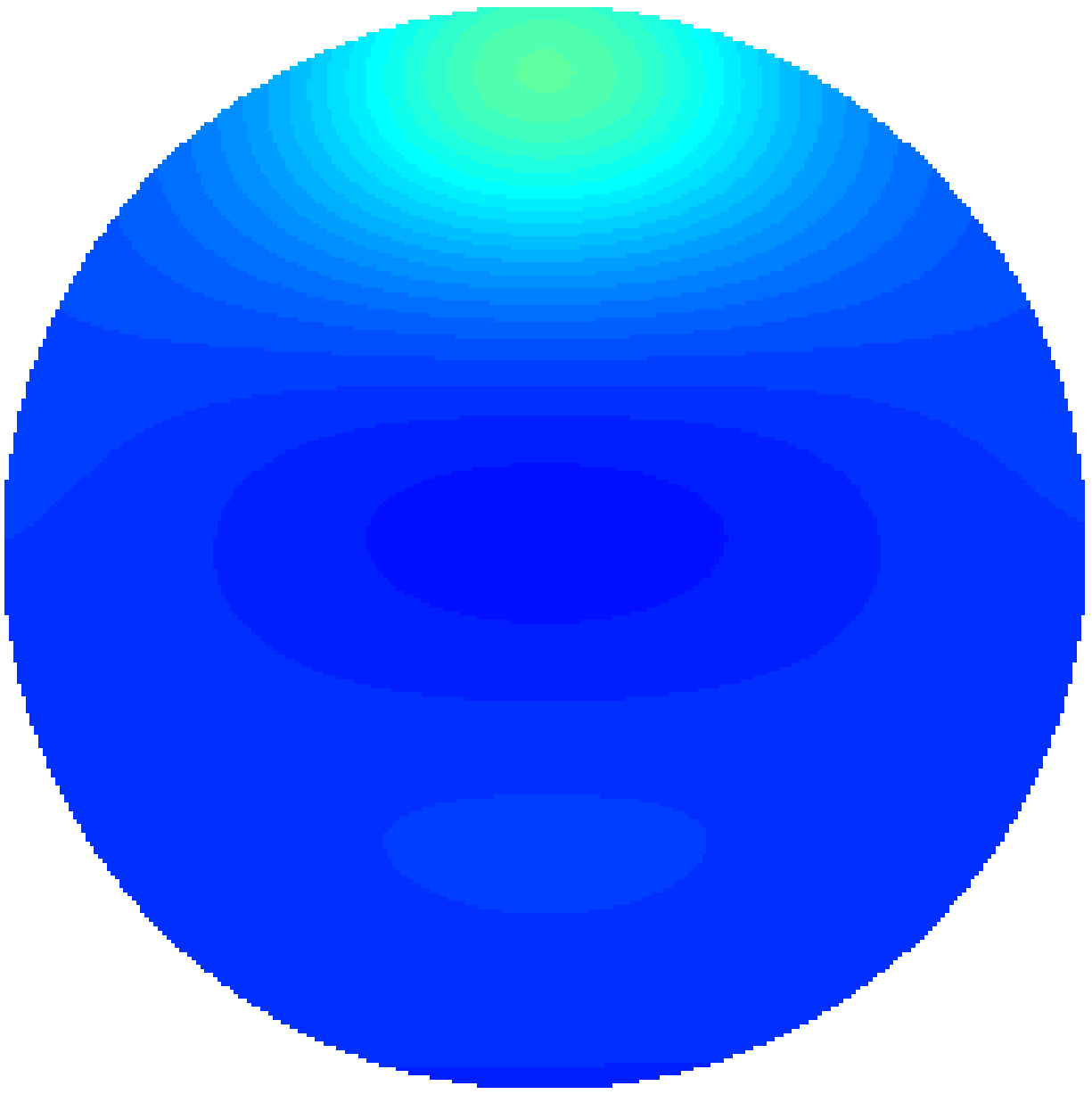}}
\epsfxsize=3cm
\put(60,0){\epsffile{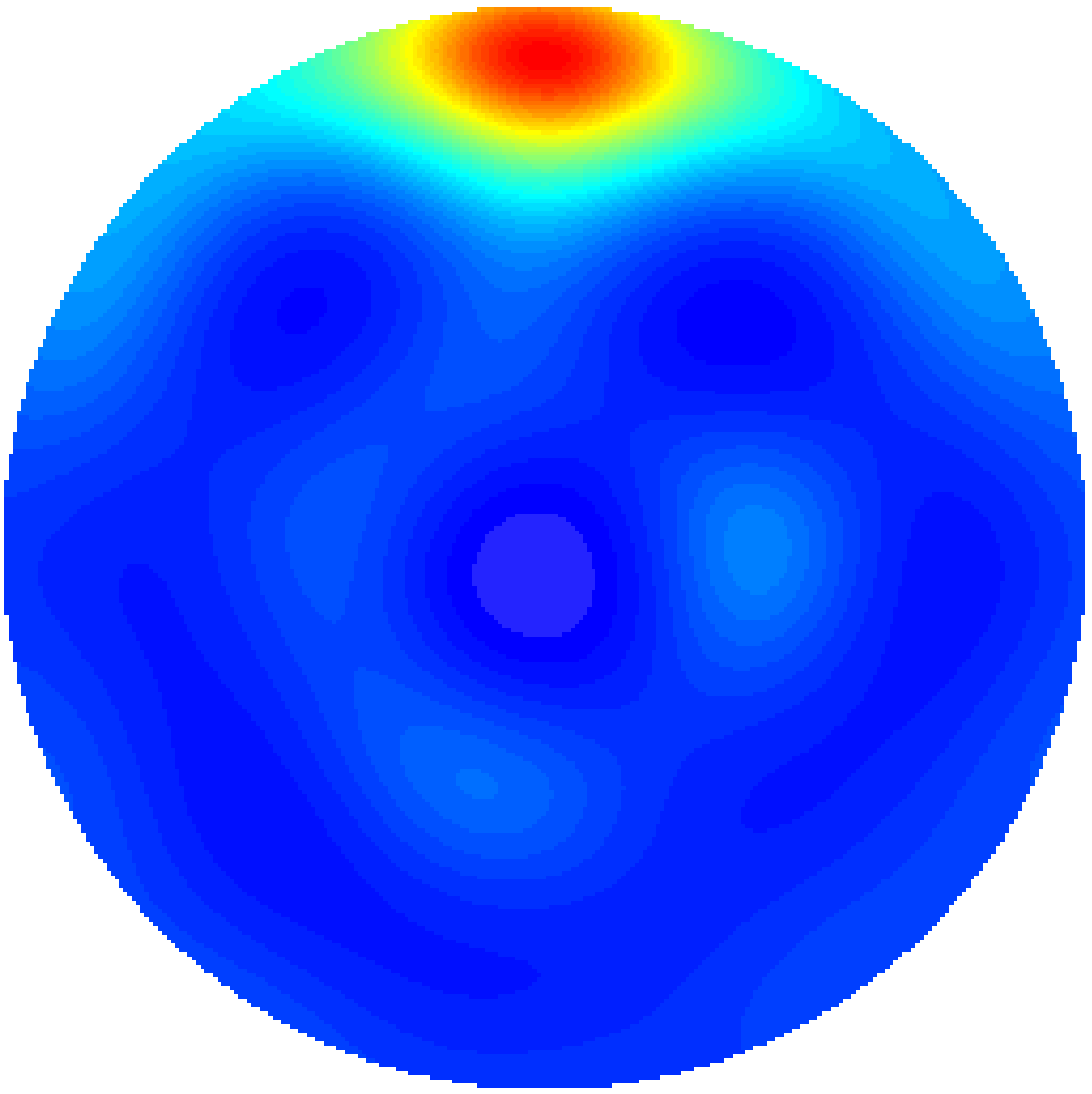}}
\epsfxsize=3cm
\put(90,0){\epsffile{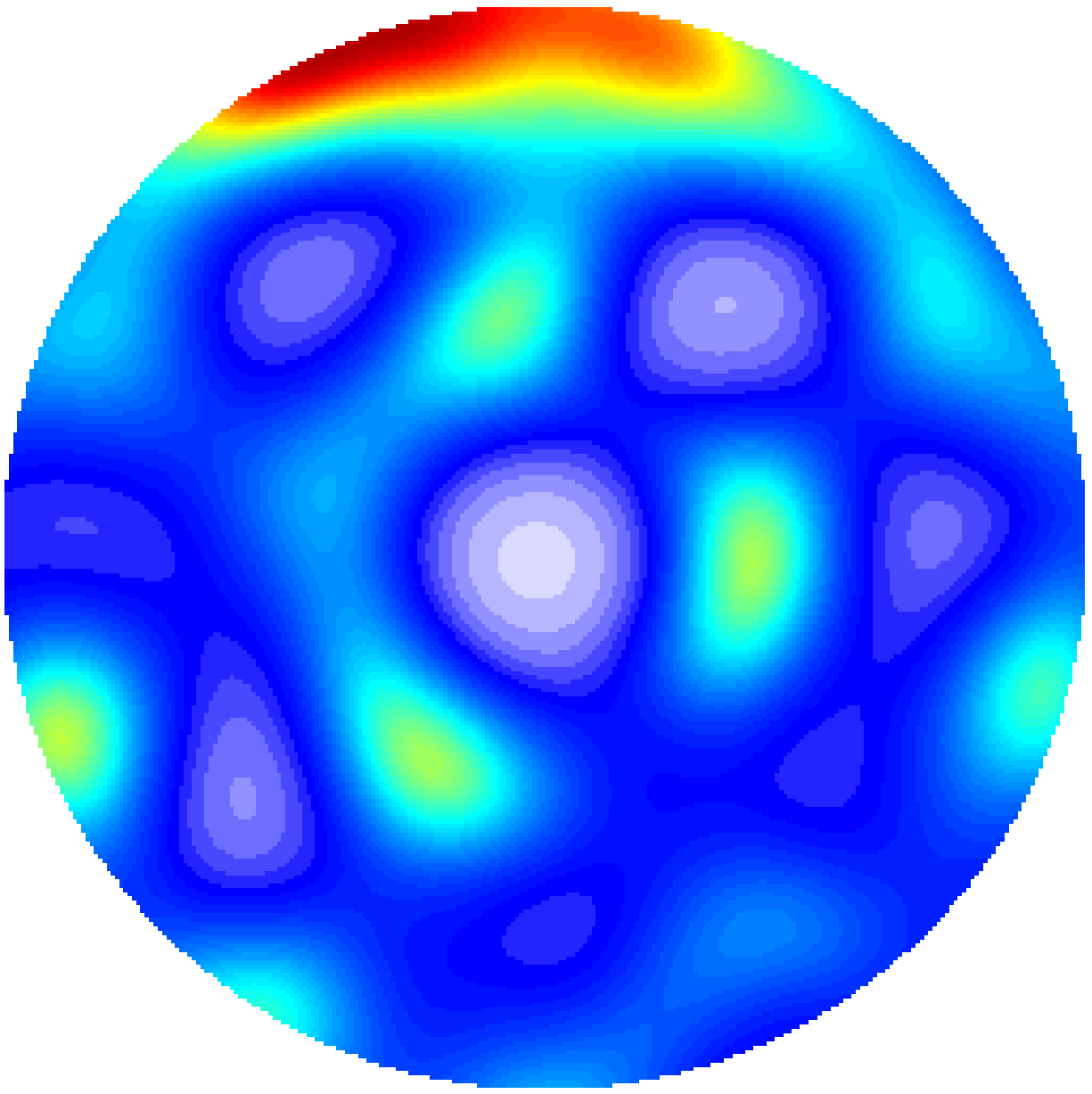}}

\put(6,59){\small Example 2}
\put(-1,55){\small Original conductivity}

\put(37,83){\small $R=3.0$}
\put(35,79){\small Uncorrected}
\put(67,83){\small $R=5.4$}
\put(65,79){\small Uncorrected}
\put(97,83){\small $R=6.0$}
\put(95,79){\small Uncorrected}

\put(37,35){\small $R=3.0$}
\put(37,31){\small Corrected}
\put(67,35){\small $R=5.4$}
\put(67,31){\small Corrected}
\put(97,35){\small $R=6.0$}
\put(97,31){\small Corrected}

\put(30,47){42\%}
\put(60,47){29\%}
\put(90,47){35\%}
\put(30,-1){35\%}
\put(60,-1){15\%}
\put(90,-1){39\%}

\end{picture}
\caption{\label{fig:ex2recon}Example 2 reconstructions; the original
  conductivity on the left, traditional D-bar reconstructions on the upper
  row and boundary corrected reconstructions on the lower row; the numbers
  beside the pictures are $L^2$ -errors, for the full error graph, see figure
  \ref{fig:errors}. The first reconstruction pair is always calculated with
$R=3$, the second one is the one with the lowest numerical $\L^2$
-error for the boundary corrected reconstruction, and the third one is with $R=6$ to show how the reconstructions fail.}
\end{figure}

\begin{figure}[p]
\begin{picture}(120,86)
\epsfxsize=3cm
\put(0,24){\epsffile{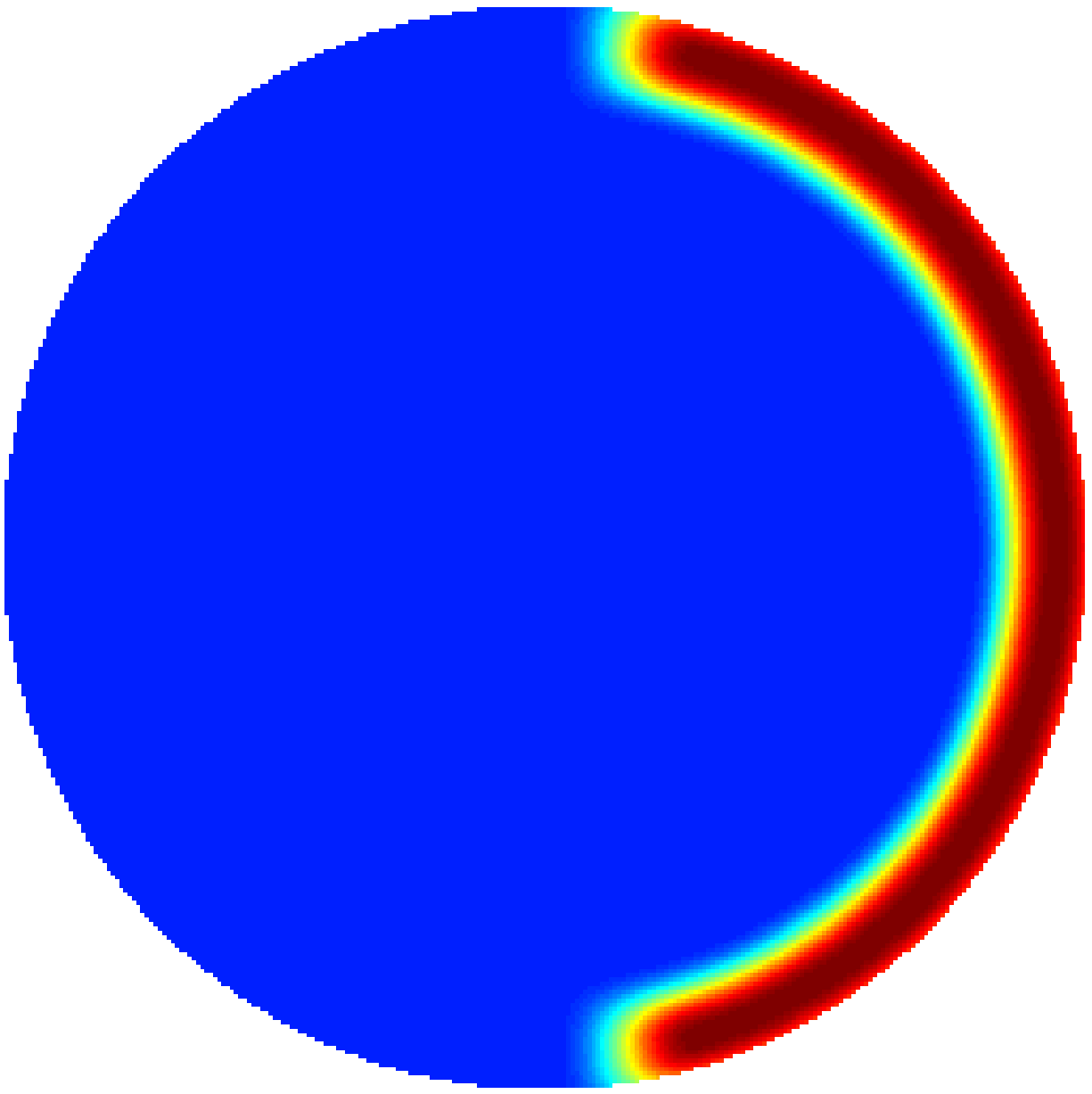}}
\epsfxsize=3cm
\put(30,48){\epsffile{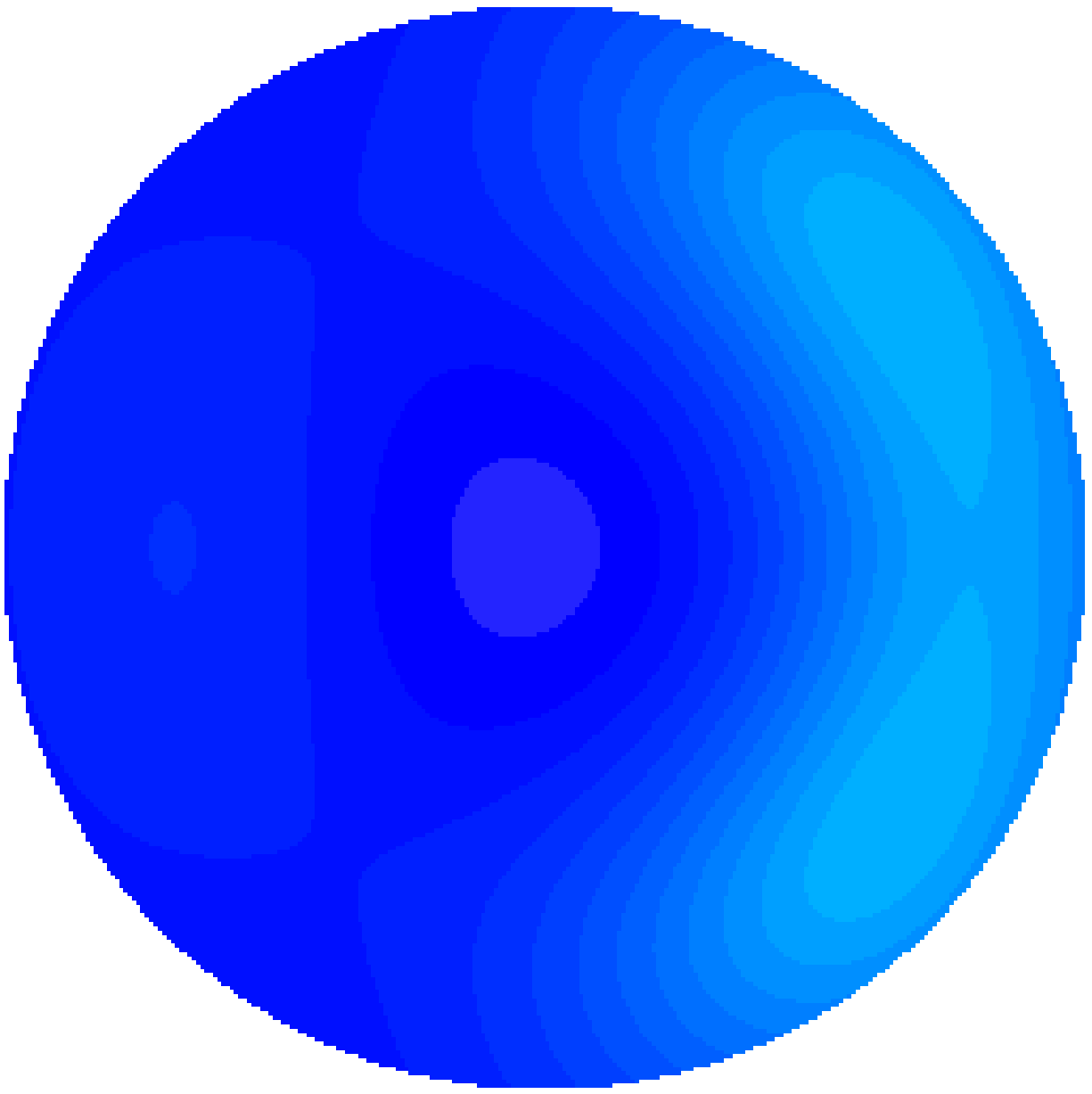}}
\epsfxsize=3cm
\put(60,48){\epsffile{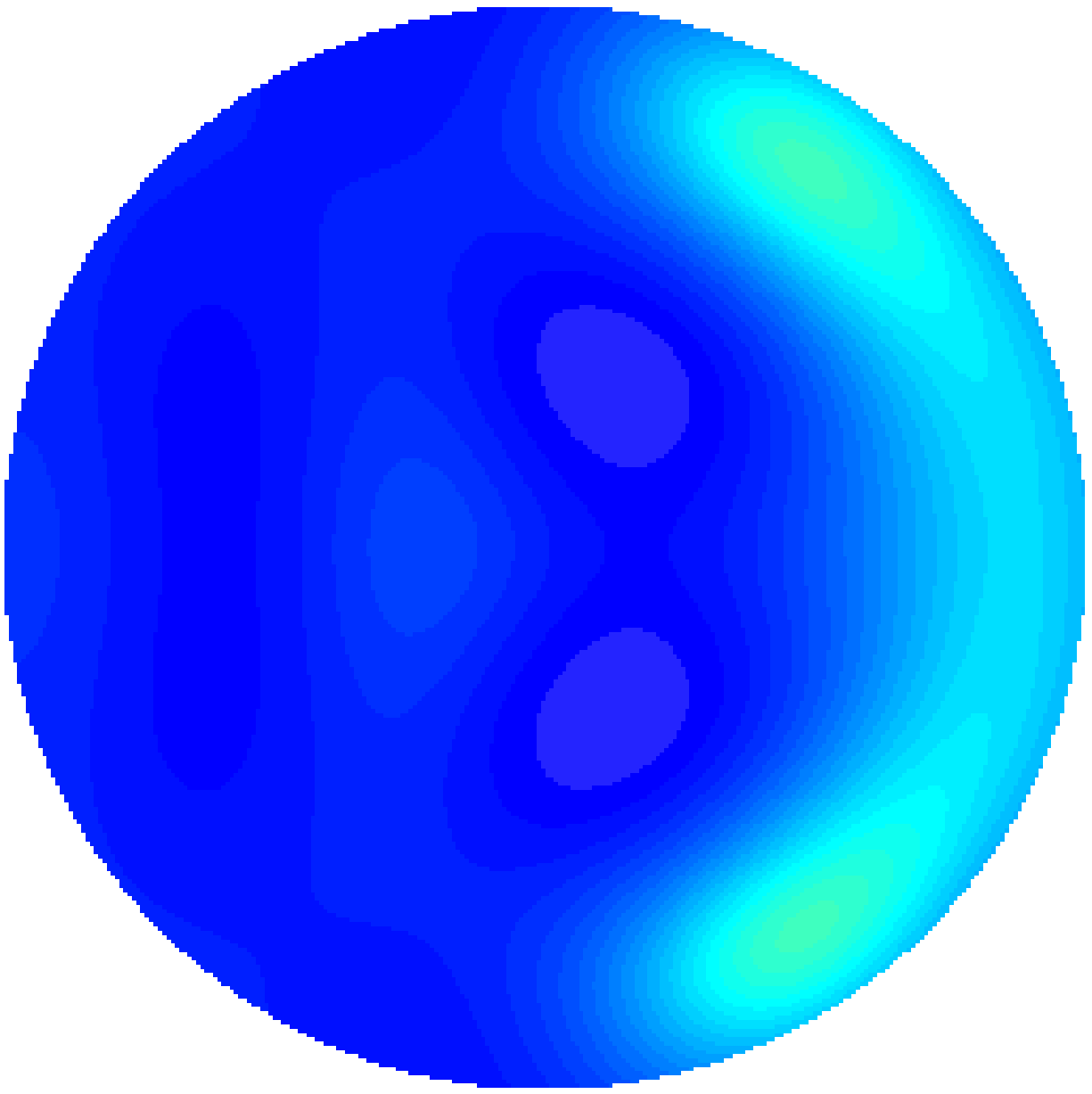}}
\epsfxsize=3cm
\put(90,48){\epsffile{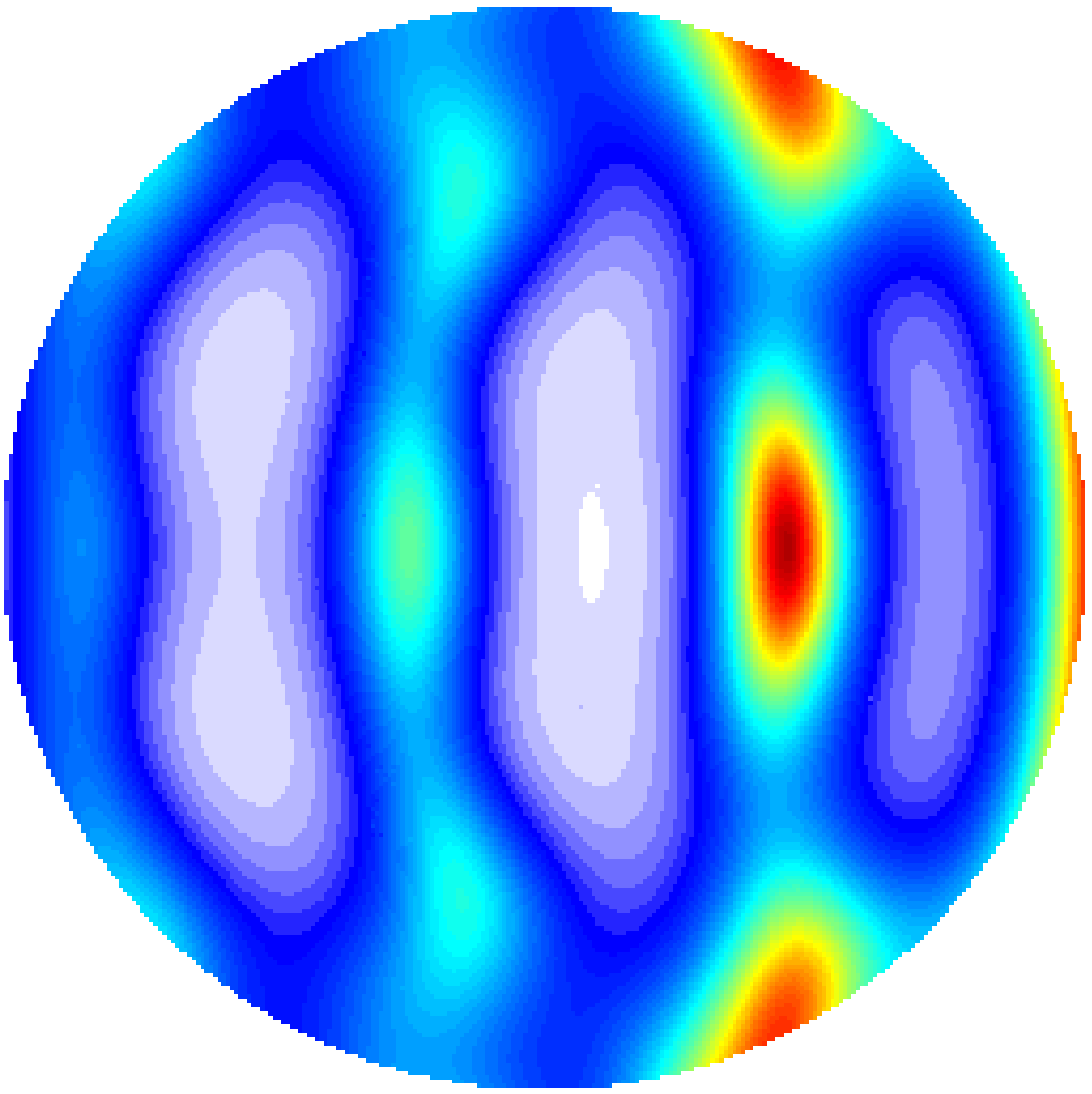}}
\epsfxsize=3cm
\put(30,0){\epsffile{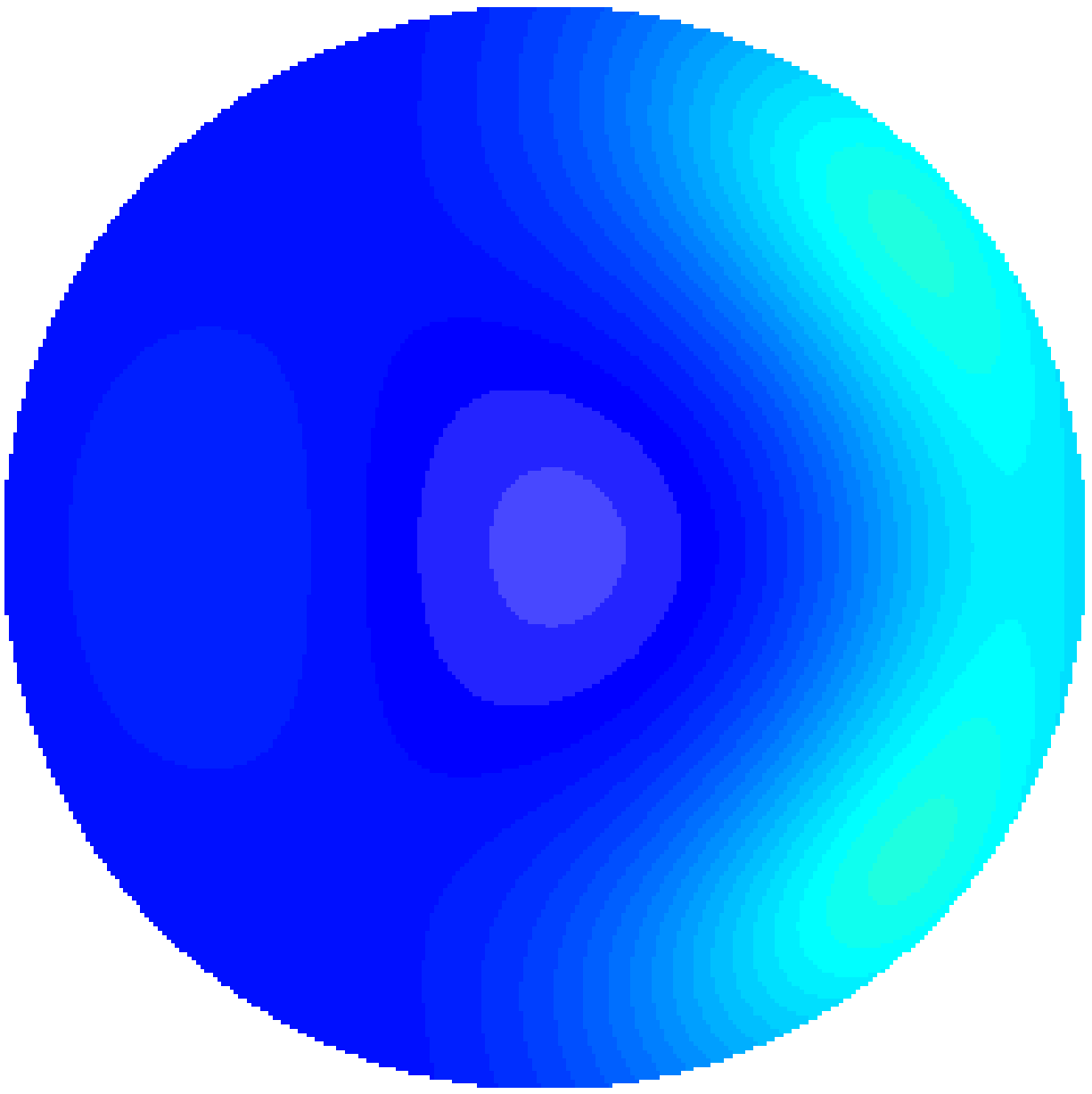}}
\epsfxsize=3cm
\put(60,0){\epsffile{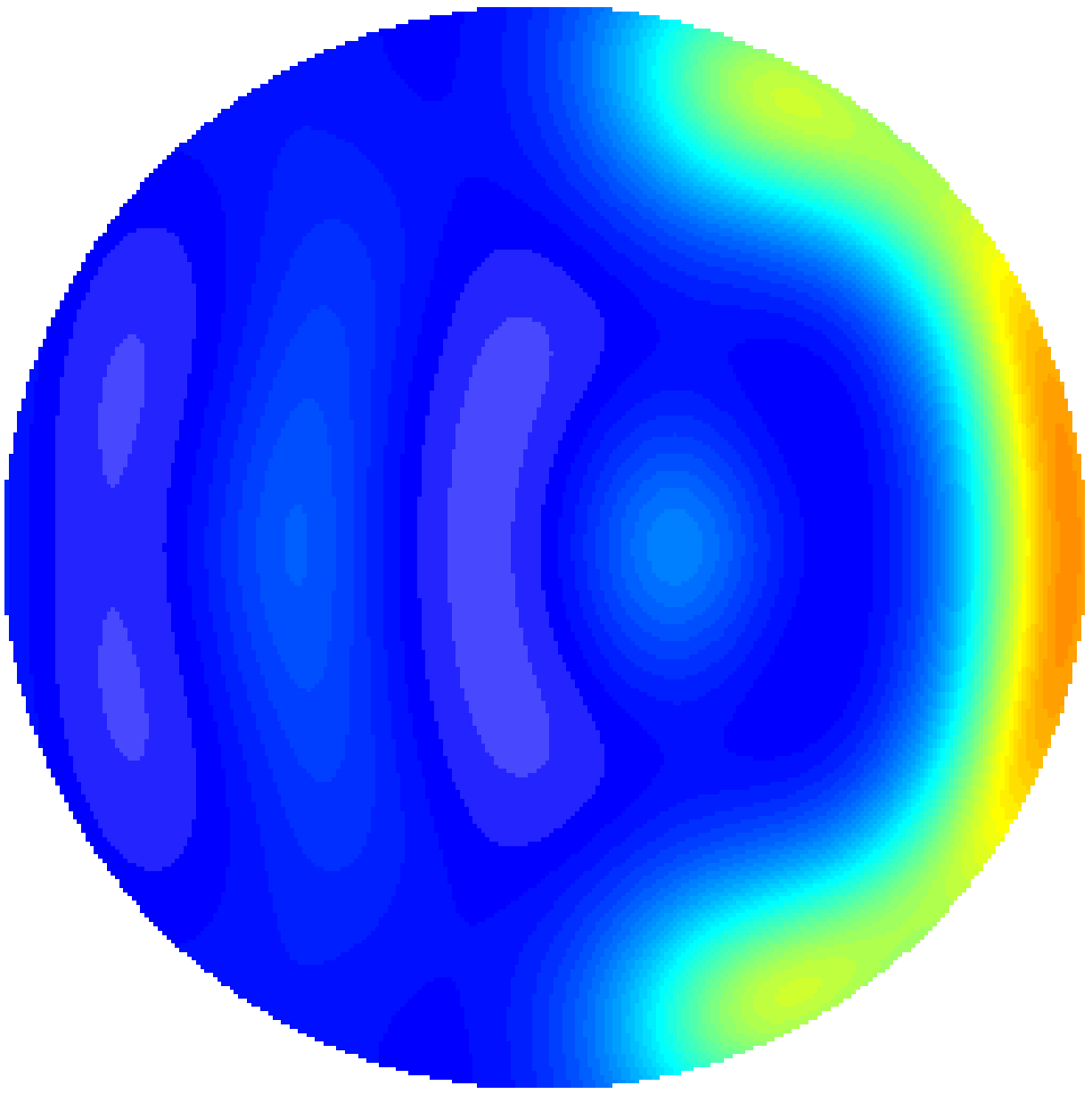}}
\epsfxsize=3cm
\put(90,0){\epsffile{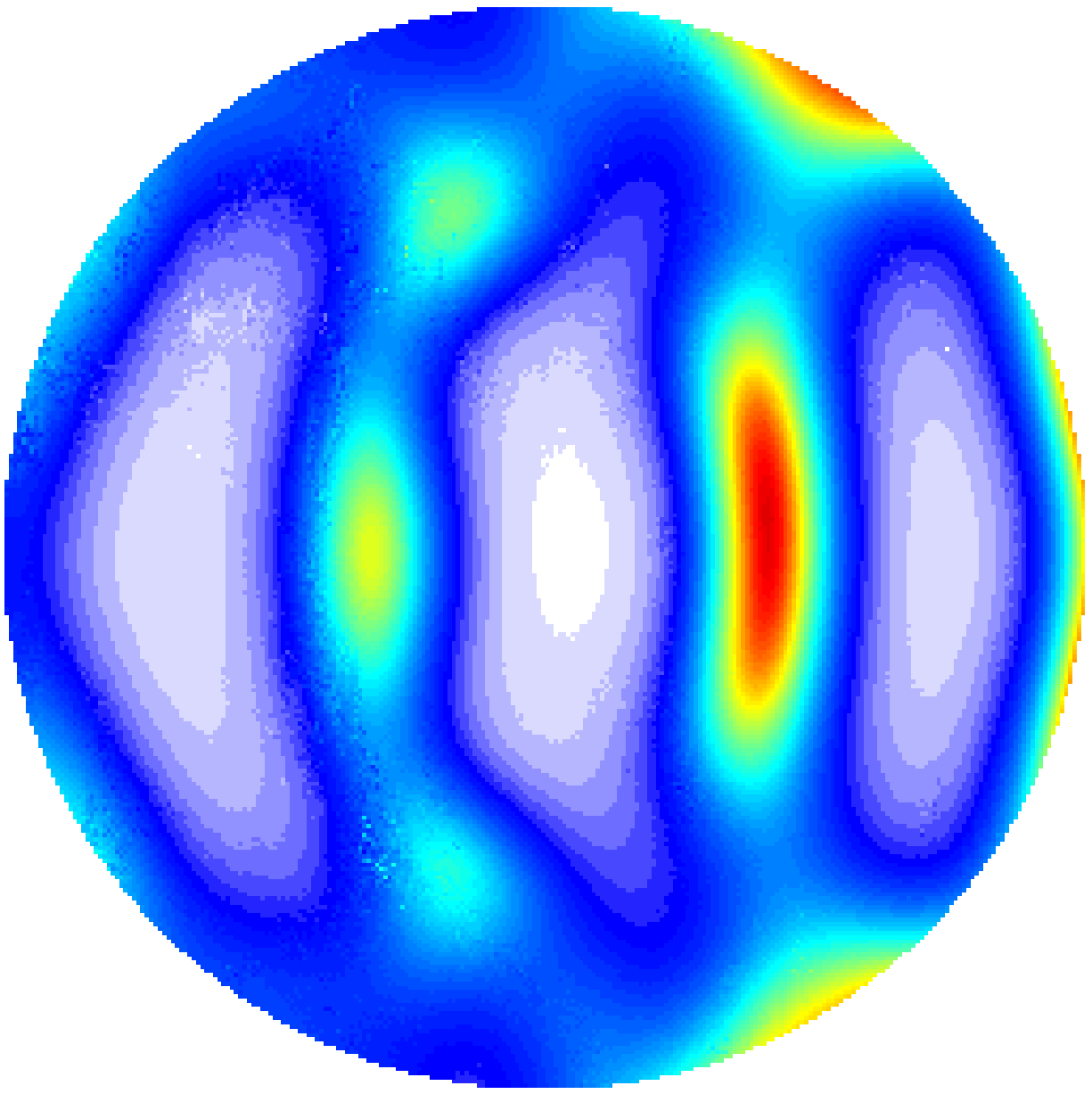}}

\put(6,59){\small Example 3}
\put(-1,55){\small Original conductivity}

\put(37,83){\small $R=3.0$}
\put(35,79){\small Uncorrected}
\put(67,83){\small $R=5.0$}
\put(65,79){\small Uncorrected}
\put(97,83){\small $R=6.0$}
\put(95,79){\small Uncorrected}

\put(37,35){\small $R=3.0$}
\put(37,31){\small Corrected}
\put(67,35){\small $R=5.0$}
\put(67,31){\small Corrected}
\put(97,35){\small $R=6.0$}
\put(97,31){\small Corrected}

\put(30,47){66\%}
\put(60,47){60\%}
\put(90,47){67\%}
\put(30,-1){59\%}
\put(60,-1){39\%}
\put(90,-1){75\%}

\end{picture}
\caption{\label{fig:ex3recon}Example 3 reconstructions; the original
  conductivity on the left, traditional D-bar reconstructions on the upper
  row and boundary corrected reconstructions on the lower row; the numbers
  beside the pictures are $L^2$ -errors, for the full error graph, see figure
  \ref{fig:errors}. The first reconstruction pair is always calculated with
$R=3$, the second one is the one with the lowest numerical $\L^2$
-error for the boundary corrected reconstruction, and the third one is with $R=6$ to show how the reconstructions fail.}
\end{figure}

\begin{figure}[p]
\begin{picture}(120,86)
\epsfxsize=3cm
\put(0,24){\epsffile{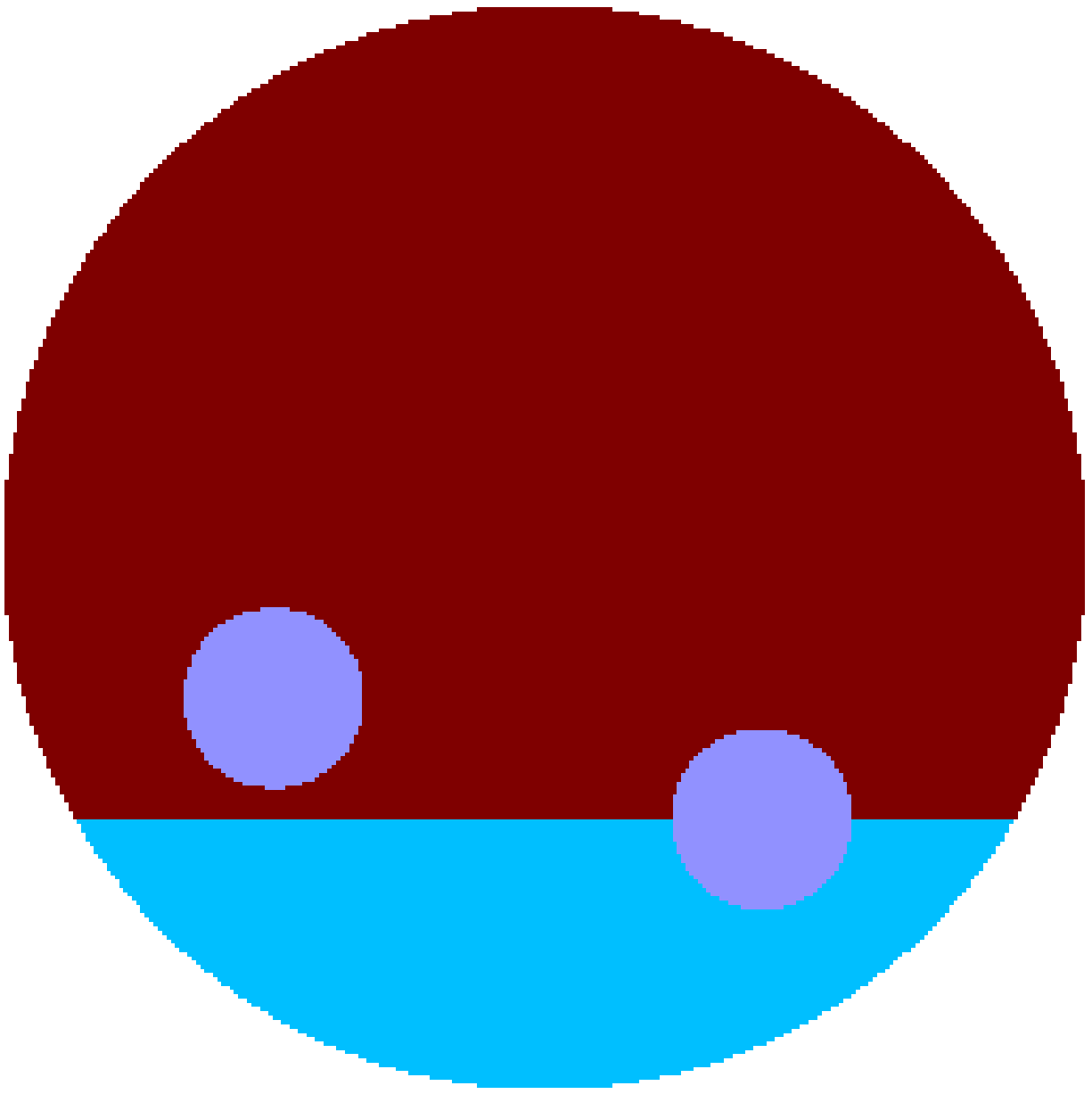}}
\epsfxsize=3cm
\put(30,48){\epsffile{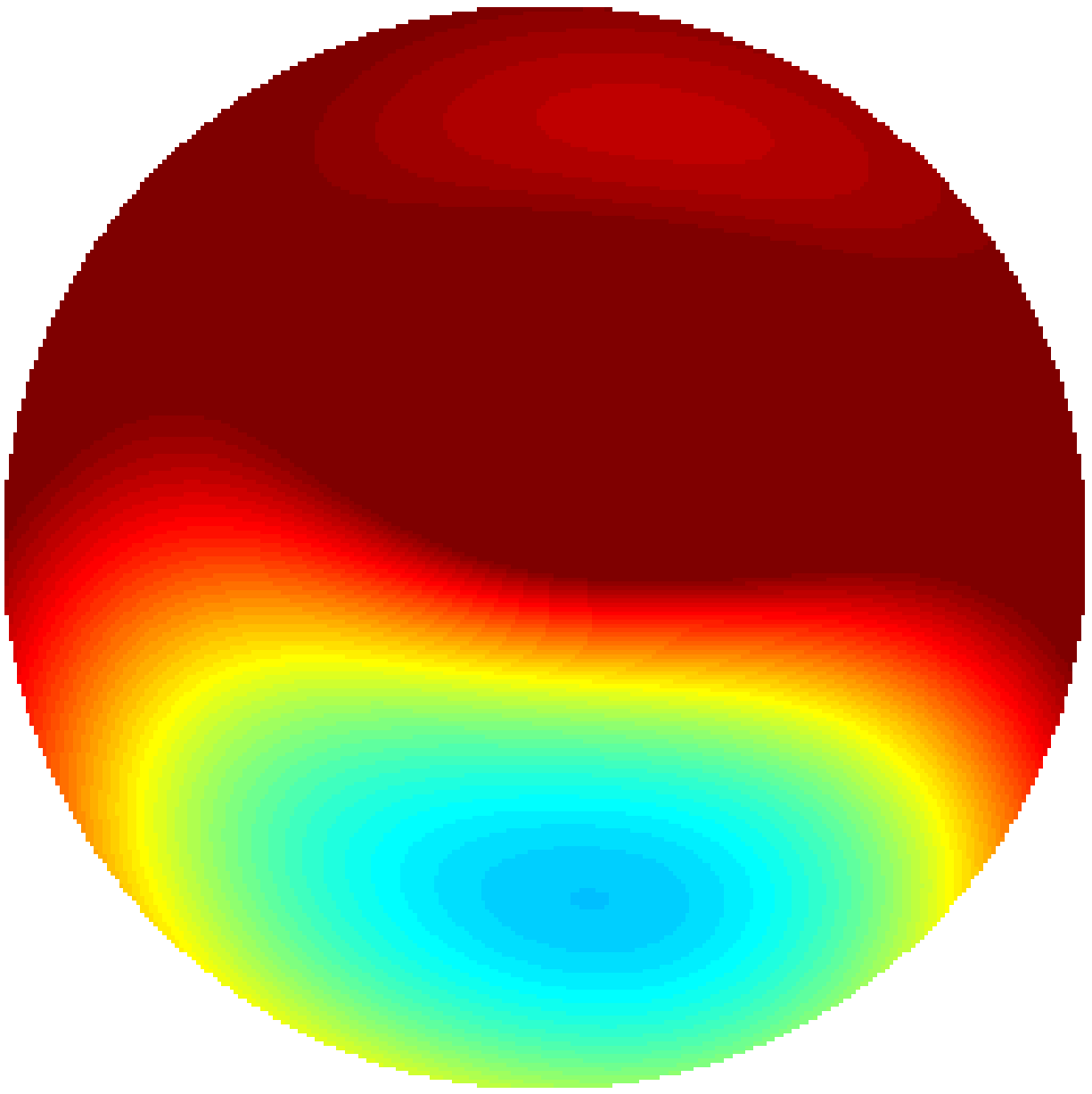}}
\epsfxsize=3cm
\put(60,48){\epsffile{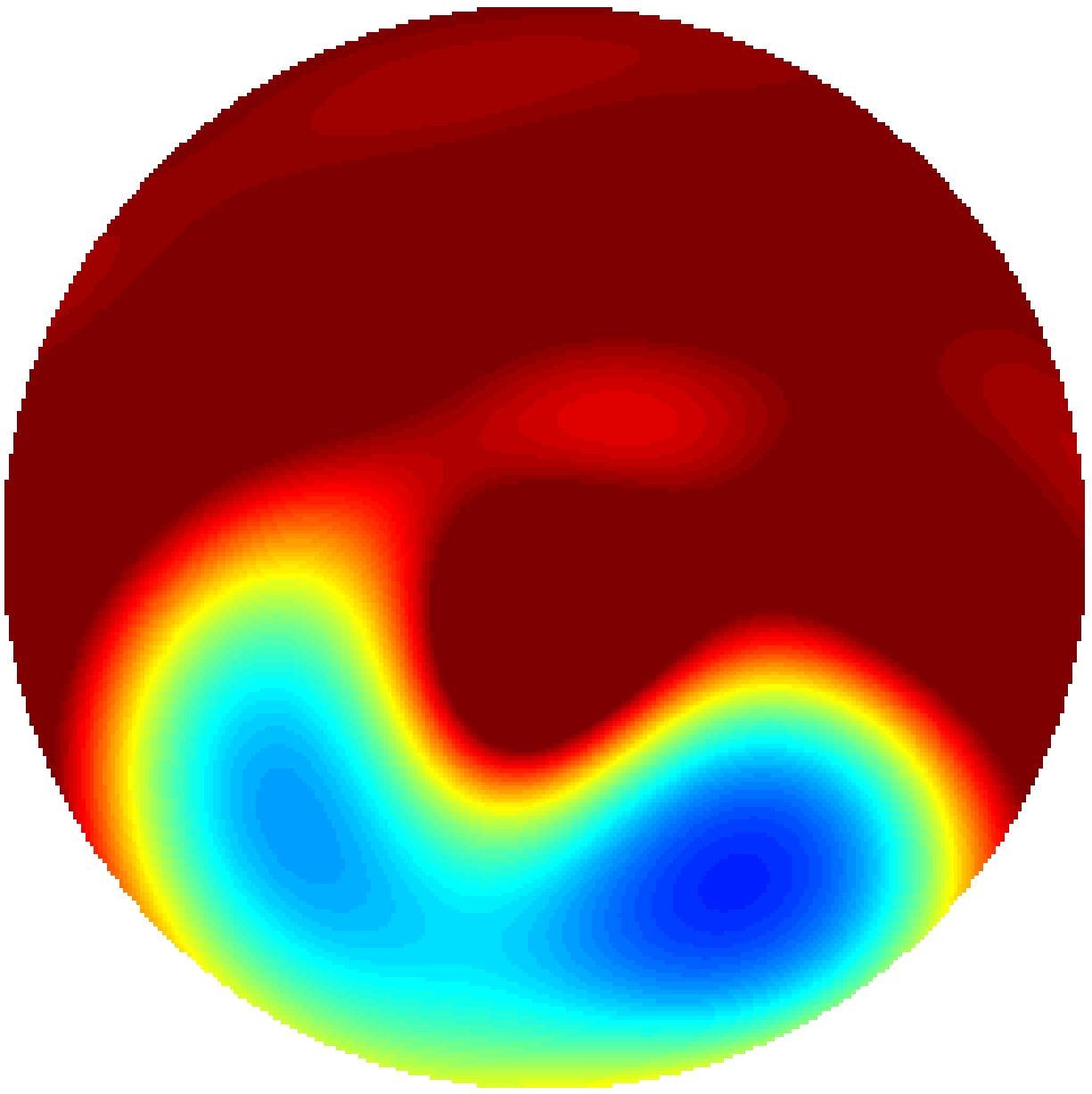}}
\epsfxsize=3cm
\put(90,48){\epsffile{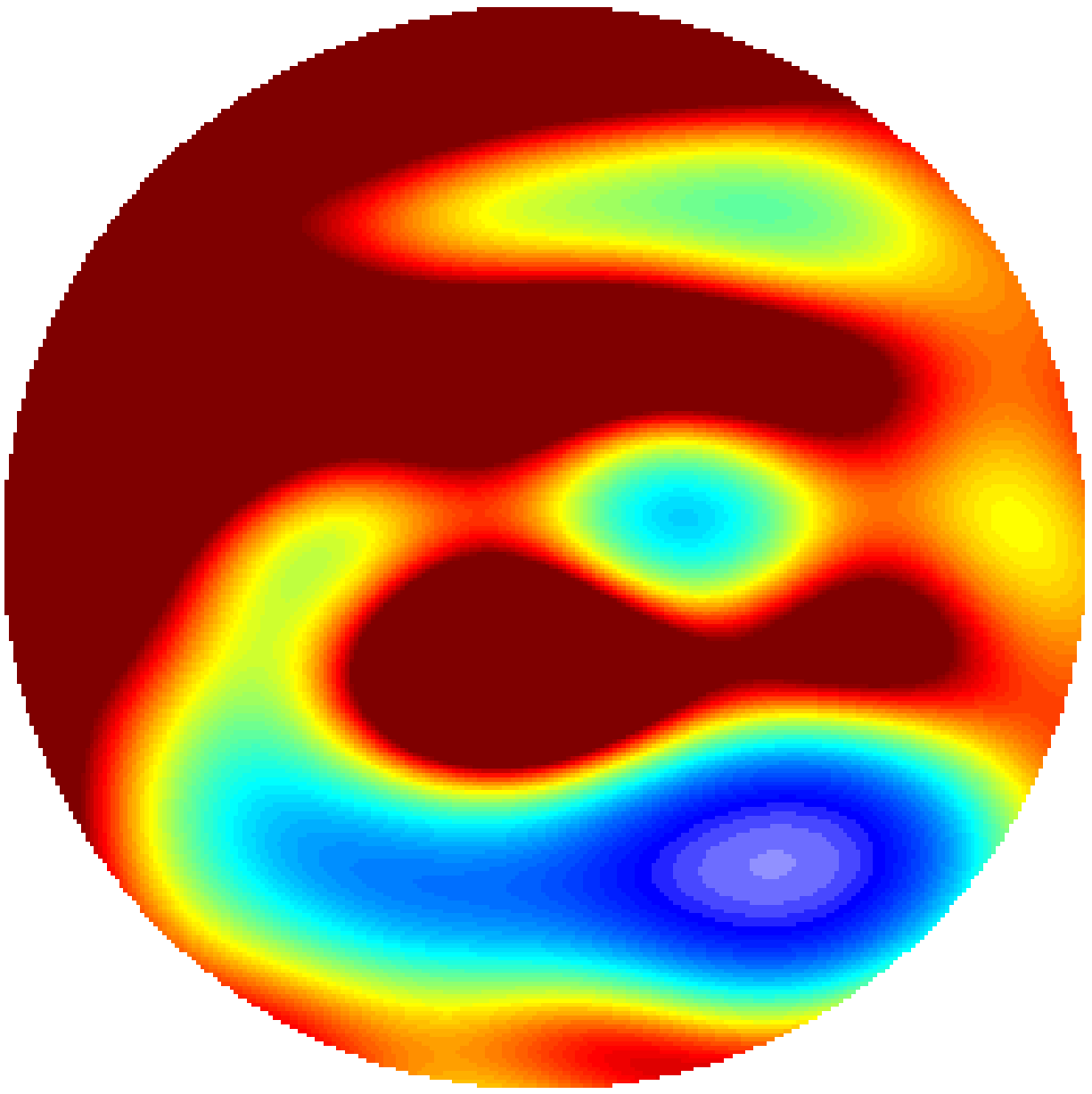}}
\epsfxsize=3cm
\put(30,0){\epsffile{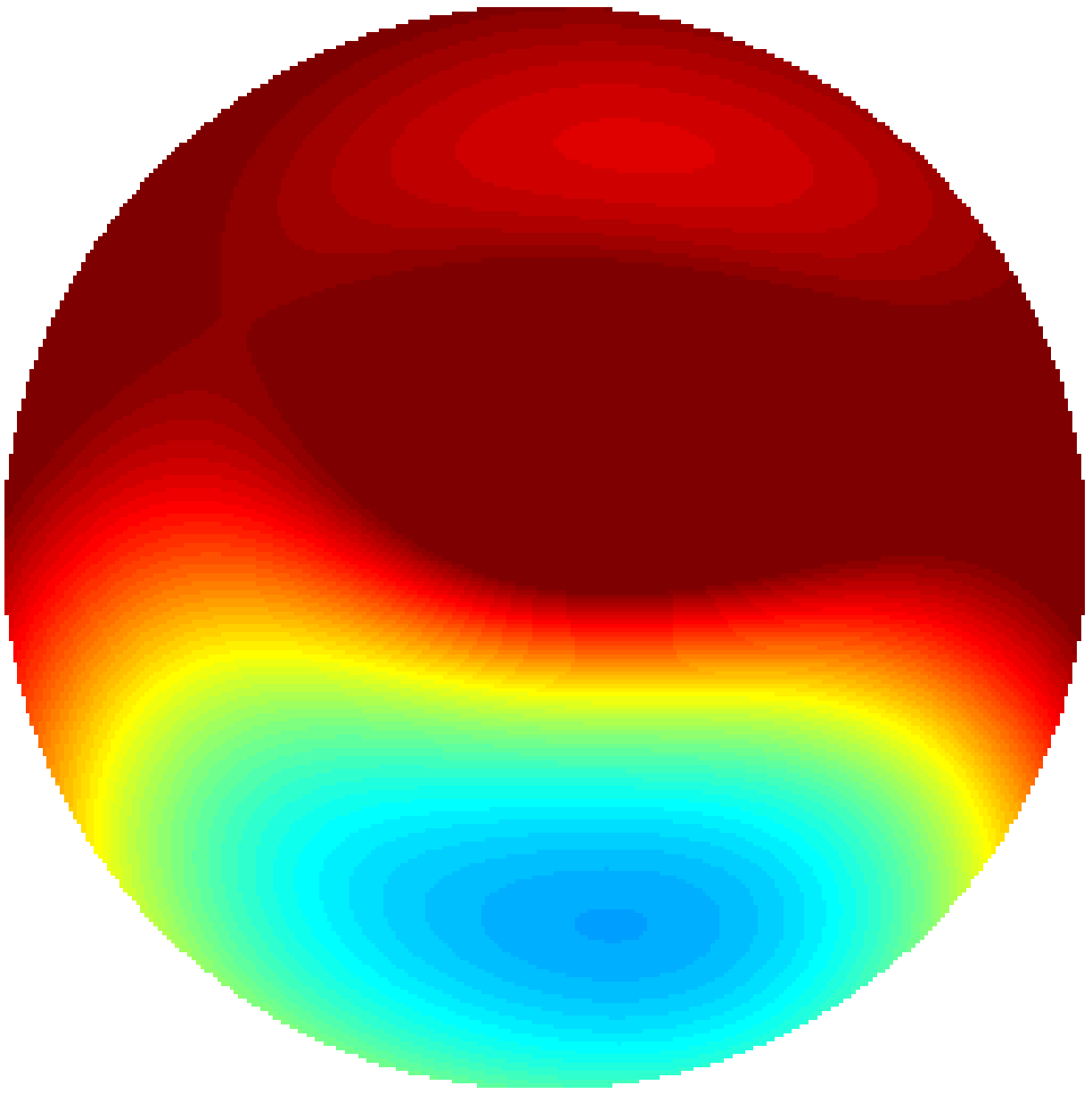}}
\epsfxsize=3cm
\put(60,0){\epsffile{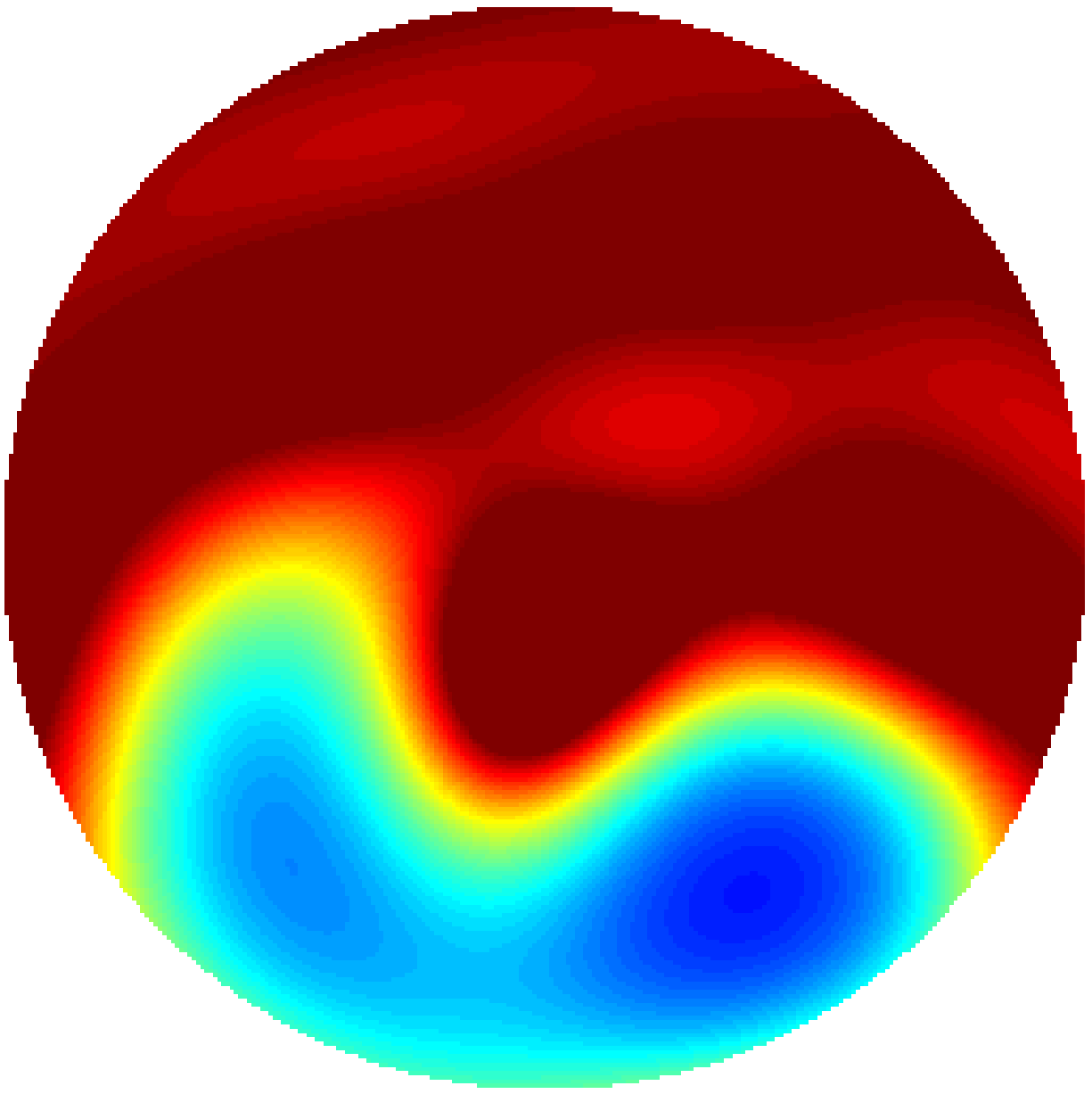}}
\epsfxsize=3cm
\put(90,0){\epsffile{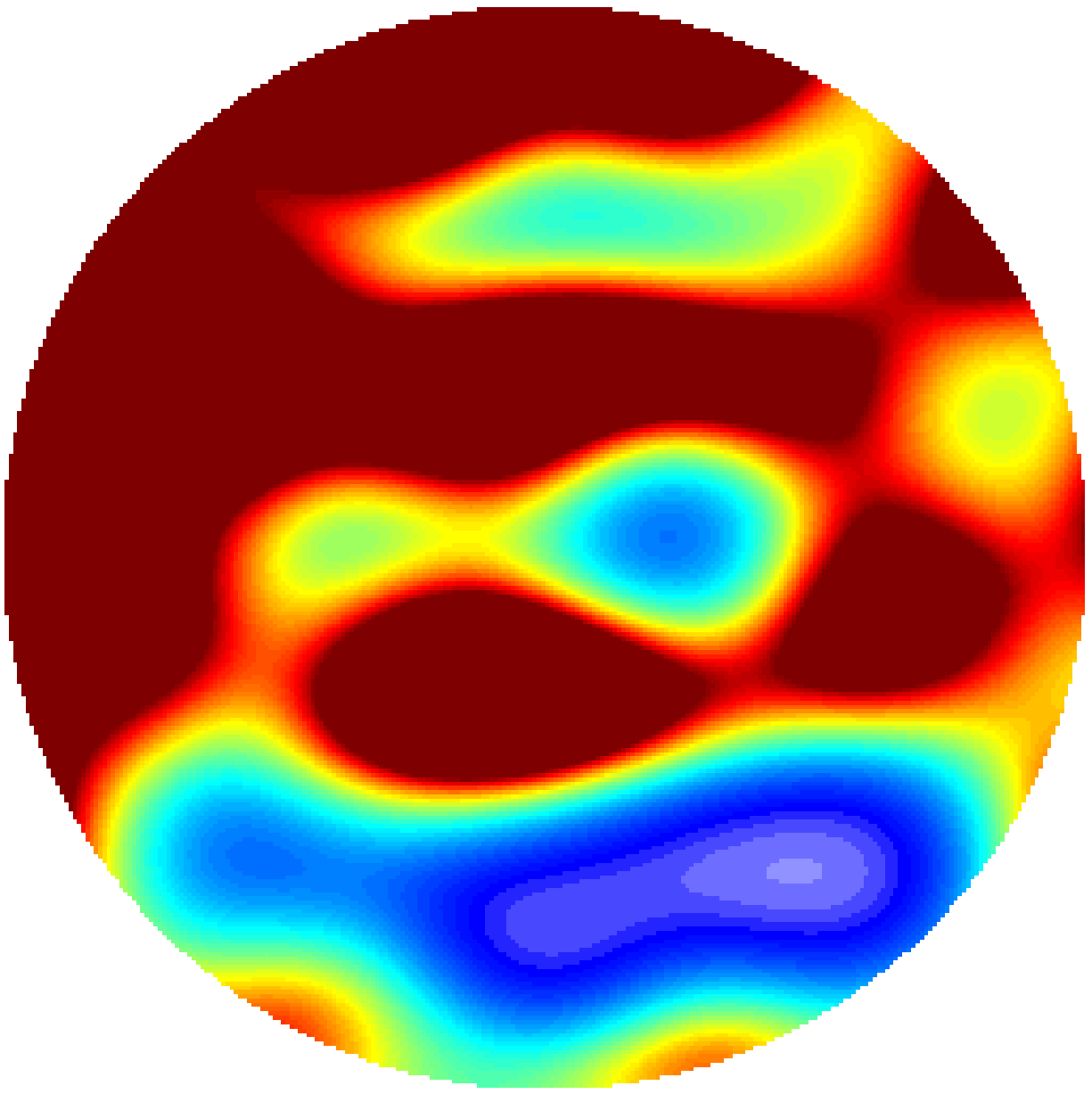}}

\put(6,59){\small Example 4}
\put(-1,55){\small Original conductivity}

\put(37,83){\small $R=3.0$}
\put(35,79){\small Uncorrected}
\put(67,83){\small $R=4.8$}
\put(65,79){\small Uncorrected}
\put(97,83){\small $R=6.0$}
\put(95,79){\small Uncorrected}

\put(37,35){\small $R=3.0$}
\put(37,31){\small Corrected}
\put(67,35){\small $R=4.8$}
\put(67,31){\small Corrected}
\put(97,35){\small $R=6.0$}
\put(97,31){\small Corrected}

\put(30,47){25\%}
\put(60,47){22\%}
\put(90,47){49\%}
\put(30,-1){25\%}
\put(60,-1){21\%}
\put(90,-1){63\%}
\end{picture}
\caption{\label{fig:ex4recon}Example 4 reconstructions; the original
  conductivity on the left, traditional D-bar reconstructions on the upper
  row and boundary corrected reconstructions on the lower row; the numbers
  beside the pictures are $L^2$ -errors, for the full error graph, see figure
  \ref{fig:errors}. The first reconstruction pair is always calculated with
$R=3$, the second one is the one with the lowest numerical $\L^2$
-error for the boundary corrected reconstruction, and the third one is with $R=6$ to show how the reconstructions fail.}
\end{figure}

For all truncation radii $R=3.0 , 3.2,\ldots ,5.8 , 6.0$, as
explained in Section \ref{sec:regDbar}, a reconstruction is calculated
with and without the boundary correction procedure using the same reconstruction points. The
conductivities and their extensions are pictured in figure
\ref{fig:conductivities}. Full error graph showing $L^2$ -error for
every reconstruction is pictured in figure \ref{fig:errors}. 
Reconstructions and the corresponding errors are pictured in figures
\ref{fig:ex1recon},\ref{fig:ex2recon},\ref{fig:ex3recon} and
\ref{fig:ex4recon}. The first reconstruction pair is always calculated with
$R=3$, the second one is the one with the lowest numerical $L^2$
-error for the boundary corrected reconstruction, and the third one is with $R=6$ to show how the reconstructions fail.

\section{Conclusion}\label{sec:conclusion}

\noindent
Our aim in this work is to find examples of simulated conductivities
that (i) share features of conductivities appearing in applications of
electrical impedance tomography, and (ii) allow higher-quality
reconstruction when boundary correction step is added to the D-bar
method. After experimenting with a large number of candidate
conductivities we concluded that for conductivities which behave
moderately at and near $\bound_1$, the method of approximating the
trace of conductivity by an optimal constant is good enough. More
precisely, the errors caused by measurement noise in
Steps (a$^\prime$) and (c$^\prime$) prevented the boundary correction
procedure from enhancing the reconstructions.

However, we were able to find several examples where the boundary
corrected D-bar method does provide better imaging quality than the
non-corrected method both in terms of quantitative error and visual
inspection.  Four such examples are presented in Section
\ref{sec:numerical}, and all of them have high contrast features in
the conductivity right at the boundary. Consequently, most medical applications do not need the boundary
correction procedure, but it may be beneficial or even necessary for
some nondestructive testing, industrial process monitoring or geophysical sensing applications.

\section*{Acknowledgments} 
\noindent
During part of the preparation of this work, SS worked as professor and JT worked as an assistant at
the Department of Mathematics of Tampere University of Technology.  
The research work of both authors was funded in part by the Finnish
Centre of Excellence in Inverse Problems Research (Academy of Finland
CoE-project 213476). The authors thank Jennifer Mueller for her
valuable comments on the manuscript.  JT was supported
in part by Pirkanmaan kulttuurirahasto. 

\end{document}